\magnification=1100

\font\bb=msbm10
\font\bbsmall=msbm7
\font\got=eufm10
\font\twelvebf=cmbx10 at 12pt

\def\N{\hbox{\bb N}}

\def\R{\hbox{\bb R}}
\def\Z{\hbox{\bb Z}}

\newdimen\sumdim
\def\diagram#1\enddiagram
{\vcenter{\offinterlineskip
\def\tvi{\vrule height 10pt depth 10pt width 0pt}
\halign{&\tvi\kern 5pt\hfil$\displaystyle##$\hfil\kern 5pt
\crcr #1\crcr}}}

\def\hfl[#1][#2][#3]#4#5{\kern-#1
\sumdim=#2 \advance\sumdim by #1 \advance\sumdim by #3
\smash{\mathop{\hbox to \sumdim{\rightarrowfill}}
\limits^{\displaystyle#4}_{\displaystyle#5}}
\kern-#3}

\def\vfl[#1][#2][#3]#4#5%
{\sumdim=#1 \advance\sumdim by #2 \advance\sumdim by #3
\setbox1=\hbox{$\left\downarrow\vbox to .5\sumdim{}\right.$}
\setbox1=\hbox{\llap{$#4$}\box1
\rlap{$#5$}}
\vcenter{\kern-#1\box1\kern-#3}}

\def\vspace[#1]{\noalign\vskip #1}

\centerline{\twelvebf Bounds on primitives of differential forms and group cocycles}
\bigskip
\centerline{Jean-Claude Sikorav}
\medskip\centerline{August 24, 2004}
\vskip10mm
\noindent{\bf 0. Introduction}
\medskip
We investigate the relations between the existence of a ``primitive'' with given bounds and the satisfaction of weighted isoperimetric inequalities. In one direction, the relation follows from various versions of Stokes' formula. In the other, it uses a Hahn-Banach type argument.
We shall consider three frameworks:
\medskip
1) Riemannian manifolds and primitives of exact differential forms. 
If $V$ is a Riemannian manifold and $\omega\in{\cal E}^q(V)$ is a differential form   of degree $q$, we define its  norm at the point $x\in V$ of by $||\omega||(x)=\max\{\omega_x(v_1,\cdots,v_q)\mid v_i\in T_xV, ||v_i||\le1\}.$

\medskip
{\bf Question 1.}\hskip3mm {\it Let $V$ be a Riemannian manifold. 
Let
$\omega\in{\cal E}^q(V)$ be exact, of degree $q\ge2$ and let $\varphi:V\to\R_+$ be a continuous function. 
When does there exist $\tau\in{\cal E}^{q-1}(V)$ such that $d\tau=\omega$ and $||\tau||\le\varphi$ ?}
\smallskip
A case of special interest will be $V=\widetilde M$, the universal covering of a compact Riemannian manifold, and $\omega$ comes from a closed form on $M$.
\medskip
2) Cellular [in particular simplicial] complexes and primitives of exact cochains.
\medskip

{\bf Question 2.}\hskip3mm {\it Let $X$ be a cellular complex. Let $u\in C^q(X;\R)$ be an exact $q$-cochain for some $q\ge2$, and let $f\in C^{q-1}(X;\R_+)$ be a nonnegative cellular $(q-1)$-cochain (function on the $(q-1)$-cells). When does there exist $t\in C^{q-1}(X;\R)$ such that $dt=u$ and $|t|\le f$ ?}
\medskip
The answer to Question 2 is an immediate application of Hahn-Banach.
\medskip
3) Groups and primitives of exact cochains.

\medskip
{\bf Question 3.}\hskip3mm {\it Let $(G,S)$ be a group equipped with a finite generating system. Let  $b$ be a $q$-cocycle on $G$ for some $q\ge2$, and let $F$ be a function from $G$ to $\R_+$.
When does there exist a $(q-1)$-cochain $a\in C^{q-1}(G;\R)$ such that $da=b$ and 
$$|a(g,g\overline s_1,g\overline s_1\overline s_2,\cdots,g\overline s_1\cdots \overline s_{q-1})|\le F(g) ?$$}
\medskip
{\bf Special case $q=2$.}\hskip3mm {\it Let $b:G^3\to\R$ be a $2$-cocycle, ie
$b(g_1,g_2,g_3)-b(g_0,g_2,g_3)+b(g_0,g_1,g_3)-b(g_0,g_1,g_2)=0$, and let $F$ be a nonnegative function on $G$. When does there exist $a:G^2\to\R$ such that $a(g_1,g_2)-a(g_0,g_2)+a(g_0,g_1)=b(g_0,g_1,g_2)$ and 
$$|a(g,g\overline s)|\le F(g) ?$$}

\bigskip
We first answer Question 1 in terms of weighted isoperimetric inequalities given by Stokes' formula. 
\medskip
{\bf Theorem 1.}\hskip3mm {\it Let $\omega\in{\cal E}^q(V)$ with $q\ge2$, and let $\varphi:V\to\R_+$ be continuous. Assume that
for every real smooth singular $q$-chain $c$ one has
$$I_c(\omega)\le{\rm M}_\varphi(I_{\partial c}).$$
Then for every $\varepsilon>0$, there exists $\tau\in{\cal E}^{q-1}(V)$ such that $d\tau=\omega$ and $||\tau||\le\varphi+\varepsilon$.
}

\medskip
Here $I_c$ is the integration current associated with $c$, and ${\rm M}_\varphi(T)$ its {\it weighted mass} of a current (see the definitions in section 1). For instance, if $\varphi=1$ and $\partial c$ has no geometric cancellations, ${\rm M}_\varphi(I_{\partial c})$ is its $(q-1)$-dimensional volume.
\medskip
{\bf Corollary.}\hskip3mm {\it The ``smallest'' norm of a primitive of $\omega$ is
$$\inf\{||\tau||_\infty\mid d\tau=\omega\}=\sup_T{T(\omega)\over{\rm M}(\partial T)}
=\sup_c{I_c(\omega)\over{\rm M}(I(\partial c)}.$$}
 \medskip
In the case of volume forms, the result is much nicer.
\medskip

{\bf Theorem 2.}\hskip3mm{\it Let $V$ be an oriented Riemannian manifold of dimension $n$
and let $\omega$ be a nonnegative smooth $n$-form (in particular a volume form).
Let $\varphi=V\to\R_+$ be continuous. 
Assume that, for every compact domain $\Omega\subset V$ with smooth 
boundary, 
 $$\int_\Omega\omega\le{\rm vol}_\varphi(\partial\Omega)= \int_{\partial \Omega}\ \varphi d\sigma$$
where $d\sigma$ is the $(n-1)$-dimensional measure on $\partial\Omega$.
 
Then for every continuous $\varepsilon>0$, there exists $\tau\in{\cal E}^{n-1}(V)$ such that $d\tau=\omega$
and $||\tau||\le \varphi+\varepsilon$.
}
\medskip
 From Theorem 1 we deduce a comparison predicted by Gromov [G2, p.98] between the cofilling function and (a suitable version of) the filling area. For the definitions, see section 4.

  \medskip
{\bf Theorem 3.}\hskip3mm {\it Let $V$ be a Riemannian manifold such  that $H_1(V;\R)=0$ and  $V$ is {\rm quasihomogeneous}: there exists $C>0$ and for every
$x,y\in V$ a $C$-bilipschitz homeomorphism $h:V\to V$ with $d(h(x),y)\le C$.
\smallskip
Then
$${\rm Cof}(R)\sim{\R FA(R)\over R},$$
where $\varphi(R)\sim g(R)$ means that there exists $C_1,C_2>0$ such that $C_1\varphi(C_1R)\le g(R)\le C_2\varphi(C_2R)$.} 
\medskip
{\bf Remark.}\hskip3mm [G2] states that ${\rm Cof}(R)\sim{\rm FA}(R)/R$. The equivalence between ${\rm FA}$ and $\R FA$ for $V$ the universal covering of a compact manifold is an old question [?]. 

\medskip
Question 3 can be answered using Hahn-Banach.
 We give first the case $q=2$:   
\smallskip
{\bf Theorem 4.}\hskip3mm {\it Let $b$ be a $2$-cocycle on $G$, and let $F$ be a function from $G$ to $\R_+$. 
Then the following are equivalent:
\smallskip
(i) There exists $a\in C^1(G;\R)$ 
such that $da=b$ and
$ |a(g,g\overline s^{\pm1})|\le F(g).$
\smallskip
(ii) For every $g\in G$ and every relation $w=s_1^{\varepsilon_1}\cdots s_n^{\varepsilon_n}\in R$, one has, setting $g_i=g\overline s_1^{\varepsilon_1}\cdots\overline s_{i-1}^{\varepsilon_{i-1}}$ ($g_0=1$):
$$\big|\displaystyle\displaystyle\sum_{i=1}^nb(1,g_i,g_{i+1})-\sum_{\varepsilon_i=-1} b(1,g_{i+1},g_i)\big|
\le\sum_{i=1}^nF(g_i).$$
}

In particular, if $b$ is bounded, it has a primitive satisfying $|a(g,g\overline s^{\pm1})|\le\delta_{\hbox{\bbsmall R}}^{ab}(|g|)$, where $\delta_{\hbox{\bbsmall R}}^{ab}$
is the Abelianized and regularized Dehn function.
\vfill\eject
\noindent{\bf Contents}
\medskip
{\bf 1. Spaces of currents. Locally flat forms}
\medskip
{\bf 2. Answer to Question 1}
\medskip
{\bf 3. The case of volume forms}
\medskip
{\bf 4. Filling and cofilling invariants}
\medskip
{\bf 5. Primitives of cocycles  of degree $2$ on a group}
\medskip
{\bf 6. Relation with the $\ell_1$-norm of Gersten and the homological Dehn function}
\medskip
{\bf 7. Filling and cofilling in groups}
\medskip
{\bf 8. Relation between Questions 1 and 2}
\medskip
{\bf 9. Relation between Questions 1, 2 and 3 for $q=2$}
\vskip10mm
\noindent{\bf 1. Spaces of currents with compact support. Approximation and regularization results}
\medskip
First, 
${\cal E}(V)=\bigoplus_{q=0}^n\ {\cal E}^q(V)$ is the topological vector spaces of smooth
differential forms. Its dual space is ${\cal E}'(V)=\displaystyle\oplus\ {\cal E}'_q(V)$. Recall [dR] that ${\cal E}_q(V)$ is reflexive.
The elements of ${\cal E}'(V)$
will be called  {\it currents with compact support}. This is a slight (but usual) abuse since the topology is distinct from that induced by the space of currents (the dual of 
forms with with compact support). It will cause no preoblem since we shall never use currents without compact support.
\medskip
We now consider special subspaces of ${\cal E}(V)$.
\medskip 1) Currents of {\it finite mass}, where the mass ${\rm M}(T)$ is defined by
${\rm M}(T)=\sup\{T(\varphi)\mid ||\varphi||\le 1\}.$ By the representation theorem of Federer [F1], these are the same as compactly supported
{\it measure-type} currents:
$$T_\xi(\varphi)=\int_V\varphi_x(\xi_x)\ d\nu(x),$$
where $\xi:M\to \Lambda^q(TM)$ is a measurable field of $q$-vectors, compactly supported, and such that $||\xi||\in L^1(V)$. Note that ${\rm M}(T_\xi)=||\xi||_{L^1}$.
\medskip
{\bf Weighted mass.}\hskip3mm 
If $\varphi$ is a nonnegative function on $V$, we can define the {\it weighted mass} of $T\in {\cal E}'_q(V)$:
$${\rm M}_\varphi(T)=\sup\{T(\omega)\mid\omega\in{\cal E}^q(V)\ ,\ ||\omega||\le\varphi\}.$$
In particular, if $f=1$ this is the usual mass.
\smallskip

We denote ${\cal M}_q(V)\subset{\cal E}'_q(V)$ the subspace of measure-type currents.
Following Federer, one defines ${\rm N}(T)={\rm M}(T)+{\rm M}(\partial T)$, and calls $T$ {\it normal} if ${\rm N}(T)$ is finite. We denote by ${\cal N}_q(V)$ the space of compactly supported normal $q$-currents, and ${\cal N}_{q,K}(V)$ the space of those with support in the compact subset $K$.

\medskip
{\bf Flat chains and locally flat cochains [W] [F1] [F2].}\hskip3mm For $K\subset V$ compact, one defines the flat semi-norm
$$F_K(T)=\sup\{T(\varphi)\mid \varphi\in{\cal E}^q(V)\ , \ \max(||\varphi||_K,||d\varphi||_K)\le1\}.$$
Then, by [F1] (p. 367),
$F_K(T)=\inf\{{\rm M}(T-\partial S)+{\rm M}(S)\mid {\rm supp}(S)\subset K \}.$
One defines ${\cal F}_{q,K}(V)$ as the $F_K$-closure of ${\cal N}_{q,K}(V)$. The space of {\it  flat $q$-chains} is ${\cal F}_q(V)=\cup_K{\cal F}_{q,K}(V)$, union over all compact subsets $K$.

\smallskip
A {\it locally flat $q$-cochain} is a linear form $\ell$ on ${\cal F}_q(V)$ (or on ${\cal N}_q(V)$) which is $F_K$-bounded on every ${\cal F}_{q,K}$ (or on every ${\cal N}_{q,K}$). Such a cochain is equivalent to a  {\it locally flat form} of degree $q$, ie $\lambda\in L^{\infty}_{loc}{\cal E}^q(V)$ 
(coefficients measurable and locally 
bounded) such that there exists $\mu\in L^{\infty}_{loc}{\cal E}^{q+1}(V)$  (necessarily unique) which satisfies $T(\mu)=\partial T(\lambda)$ for every $T\in{\cal N}_k(V)$ ($d\lambda=\mu$ in the sense of distributions).

\smallskip
The correspondence $\lambda\leftrightarrow\ell$ is given by
$\ell(T_\xi)=\int_V\lambda(\xi) d\nu$ if $\xi$ is a compactly supported  field of $q$-vectors.  We define $T_\xi(\lambda)=\ell(T_\xi)$, thus $T(\lambda)$ is defined if ${\rm M}(T)<\infty$.
\smallskip We denote by ${\cal F}^*_{loc}(V)$ the space of locally flat forms.

\medskip
2) {\it Smooth currents} are currents of the form $T_\xi$ where $\xi$ is a smooth (compactly supported)  field of $q$-vectors. 
These are also called {\it diffuse} currents [Su]. 
We denote ${\cal S}'_q(V)\subset{\cal E}'_q(V)$ the subspace of smooth currents.
\medskip
3) {\it Currents associated to singular chains:} if $c=\sum_{i=1}^ka_i\sigma_i$ is a real Lipschitz singular chain, one associates the integration current
$$I_c(\varphi)=\sum_ia_i\int_{\Delta^q}\sigma_i^*\varphi.$$
Note that this time, $c\mapsto I_c$ is not injective. Note also that $I_{\partial c}={\partial I}_c$ and that
${\rm M}(I_c)\le\sum_i|a_i|{\rm vol}(\sigma_i),$
with equality if there is no geometric cancellation between the $\sigma_i$, eg if there images are disjoint.
\smallskip
We denote ${\cal C}^{\rm Lip}_q(V)\subset{\cal E}'_q(V)$ the subspace of currents associated to Lipschitz singular chains, and similarly
${\cal C}^{C^k}_q(V)$ for $k\in\N\cup\{\infty\}$.
\medskip
We shall need the following density and regularization results.
\medskip

1) {\bf Density of smooth singular chains.}\hskip3mm
{\it Let $T$ be a normal current on $V$ with support contained in the interior of a compact set $K$. Then for every $\varepsilon>0$ there exists a smooth singular $\R$-chain $c$ with values in $K$, such that
$$\eqalign{&F_K(I_c-T)<\varepsilon\cr
&{\rm M}(I_c)<{\rm M}(T)+\varepsilon\cr
&{\rm M}(\partial I_c)< {\rm M}(\partial T)+\varepsilon.\cr}$$
}
\medskip
{\bf Proof.}\hskip3mm Federer ([F1], Theorem 4.2.24) proves the 
case $V=\R^n$ and $T$ normal, with $c$ polyhedral. The last two inequalities are replaced by ${\rm N}(T')\le {\rm N}(T)+\varepsilon$, but actually he proves the more precise inequalities stated here.
\medskip
In general, we embed isometrically $i:V\to\R^N$, and work in an arbitrarily small compact tubular neighbourhood $\widehat K$ of $i(K)$, equipped with a smooth projection $\pi:\widehat K\to K$ with $||\pi-{\rm Id}||<\varepsilon,||D\pi||\le 1+\varepsilon$. Let $p$ be a polyhedral chain in $\widehat K$ satisfying $F_K(I_p-T)<\varepsilon/2$, ${\rm M}(I_p)<{\rm M}(T)+\varepsilon/2$, ${\rm M}(\partial I_p)< {\rm M}(\partial T)+\varepsilon/2$.
\smallskip
Then $I_{\pi\circ p}=\pi_*(I_p)$, ${\rm M}(\pi_*I_p)\le(1+\varepsilon){\rm M}(I_p)<{\rm M}(T)+\varepsilon$, and similarly ${\rm M}(\partial\pi_*I_p)={\rm M}(\pi_*\partial I_p)<{\rm M}(\partial T)+\varepsilon$,  $F_K(\pi_*I_p-T)=F_K(\pi_*(I_p-T))\le(1+\varepsilon)F_K(I_p-T).$
Thus $c=\pi\circ p$ is the desired singular chain.
\medskip
2) In his book [dR], de Rham proves a regularization theorem for currents. It is easy to adapt his proof (p. 72-83) in the dual setting of locally flat forms, to obtain the following result. See  also  [F2] in the case where $V$ is an open set in $\R^n$. 
\medskip
{\bf Regularization of locally flat forms.} \hskip3mm {\it Let $\rho:V\to\R_+^*$ be continuous. There exists a linear chain map of degree $0$,  
${\cal R}^*_\rho:{\cal F}^*_{loc}(V)\to{\cal E}^*(V)$, and a homotopy 
${\cal R}^*_\rho-{\rm Id}=d{\cal H}^*_\rho+{\cal H}^*_\rho d$,
with the following properties:
$$\eqalign{
||{\cal R}^*_\rho\omega(x)||&\le(1+\rho(x))\ ||\omega|B(x,\rho(x))||\cr
||{\cal H}^*_\rho\omega(x)||&\le\rho(x)\ ||\omega|B(x,\rho(x))||.\cr}$$
Also, if $\omega$ is already smooth, ${\cal R}^*_\rho\omega\to\omega$ in the $C^\infty_{loc}$ topology if $\rho\to0$ in the compact-open topology.}

\vfill\eject
\noindent{\bf 2. Answer to Question 1}
\medskip
Let $\omega\in{\cal E}^q(V)$ for some $q\ge2$. Assume that it has a primitive $\tau\in{\cal E}^{q-1}(V)$ such that
$||\tau||\le\varphi$. If $T\in{\cal E}'_q(V)$, the ``Stokes identity''
$T(\omega)=T(d\tau)=(\partial T)(\tau)$
implies the weighted isoperimetric inequality
$$T(\omega)\le {\rm M}_\varphi(\partial T).$$
In particular, if $T=I_c$ is associated to a singular chain, this is the inequality of Theorem 1. This theorem states that the converse is almost true. 
\medskip
{\bf Lemma 1.}\hskip3mm {\it Under the hypothesis of Theorem 1, 
we have
$$(\forall T\in{\cal N}_q(V))\hskip5mm T(\omega)\le {\rm M}_\varphi(\partial T).$$ }
{\bf Proof of the lemma.}\hskip3mm Let $T\in{\cal N}_q(V)$. Let $K\subset V$ be a compact set such that ${\rm supp}(T)\subset {\rm Int}(K)$. By the density of smooth singular chains, for every $\varepsilon>0$ there exists a smooth singular $\R$-chain $c$ with values in $K$, such that
$$\eqalign{&F_K(I_c-T)<\varepsilon\cr
&{\rm M}(I_c)<{\rm M}(T)+\varepsilon\cr
&{\rm M}_\varphi(\partial I_c)< {\rm M}_\varphi(\partial T)+\varepsilon.\cr}$$
The first inequality says that there exists $S$ with ${\rm M}(T-I_c-\partial S)+{\rm M}(S)<\varepsilon$. Since $\omega$ is closed,
$$\eqalign{T(\omega)&=I_c(\omega)+(T-I_c-\partial S)(\omega)\le I_c(\omega)+\varepsilon||\omega||_K\cr
&\le {\rm M}_\varphi(I_{\partial c}) + \varepsilon||\omega||_K\hskip5mm\hbox{\rm (by the hypothesis)}\cr
&\le {\rm M}_\varphi(\partial T)+\varepsilon + \varepsilon||\omega||_K.\cr}$$
Since this holds for every $\varepsilon>0$, Lemma 1 follows.

\medskip
{\bf Proof of Theorem 1.}\hskip3mm If $S\in\partial{\cal N}_q(V)$, define
$\bar t(S)=T(\omega)$ for any $T\in{\cal N}_q(V)$ such that $\partial T=S$. This is well defined since $\omega$ is exact.
Moreover, for every $S=\partial T\in\partial{\cal N}_q(V)$, Lemma 1 says that
$\bar t(S)\le  {\rm M}_\varphi(\partial T)={\rm M}_\varphi( S).$
By Hahn-Banach, $\bar t$ can be extended to a linear form $t$ on ${\cal N}_{q-1}(V)$ such that 
$$(\forall S\in{\cal N}_{q-1}(V))\hskip5mm t(S)\le {\rm M}_\varphi(S).$$
Thus $t$ is defined by a $L^\infty_{loc}$ form $\tau_0$, satisfying $||\tau_0||\le\varphi$ ae. The identity
$t(\partial T)=T(\omega)$ if $T\in{\cal N}_q(V)$ means that $\tau_0$ is locally flat and   $d\tau_0=\omega$ in the sense of distributions. 
\smallskip
Using the regularization theorem of section 1, define 
$$\tau={\cal R}^*_\rho\tau_0-{\cal H}^*_\rho\omega=\tau_0+d{\cal H}^*_\rho \tau_0$$
for some $\rho\in C^0(V,\R_+^*)$. Then $\tau$ is smooth and $d\tau=d\tau_0=\omega$. Moreover, for every $x\in V$ one has
$$||\tau(x)||\le(1+\rho(x))\big)\ ||\varphi|B(x,\rho(x))||+
\rho(x)\ ||\omega|B(x,\rho(x))||.$$
If $\rho$ decreases sufficiently fast, the right-hand side is $\le \varphi(x)+\varepsilon$ for very $x\in V$, qed.

\medskip
We now state and prove a ``localized'' generalization.\medskip
{\bf Theorem 1'.}\hskip3mm {\it Let $U\subset V$ be an open subset, and let $\omega\in{\cal E}^q(V)$ with $q\ge2$, and let $\varphi:U\to\R_+$ be continuous, where $U\subset V$ is open.
Assume that $\omega$ is exact and $I_c(\omega)\le{\rm M}_\varphi(I_{\partial c})$ for every real smooth singular $q$-chain $c$  on $V$ with boundary in $U$. 
\smallskip
Then for every $\varepsilon>0$ and every compact $A\subset U$, there exists a smooth form $\tau\in{\cal E}^{q-1}(V)$ such that $d\tau=\omega$ on $V$ and $||\tau||\le \varphi+\varepsilon$ on $A$.}
\medskip
{\bf Lemma 2.}\hskip3mm {\it Let $K\subset V$ be compact.
There exists a positive continuous function $F$ on $V$ with the following property. 

For every $q-1$-current $S_1$ on $V$ of finite mass which is homologous to a current with with support in $K$, there exists
$T_1\in{\cal N}_q(V)$ such that $\partial T_1=S_1+S_2$ with ${\rm supp}(S_2)\subset K$ and 
${\rm M}_{||\omega||}(T_1)+{\rm M}_\varphi(S_2)\le {\rm M}_F(S_1)$.}
\medskip
{\bf Proof of Theorem 1'.}\hskip3mm  Let $T$ be an element of ${\cal N}_q(V)$. 
We apply the lemma to  a compact $K\subset U$ such that $A\subset{\rm Int}(K)$, and $S_1=\partial T\setminus K$. Then
$$\eqalign{\partial (T-T_1)&=(\partial T\cap K)+S_1-\partial T_1\cr
&=(\partial T\cap K)-S_2.\cr}$$
This is supported in $K$ and a fortiori in $U$, thus by the hypothesis and Lemma 1, one has
$$\eqalign{(T-T_1)(\omega)&\le {\rm M}_\varphi((\partial T\cap A)-S_2)\cr
&\le {\rm M}_\varphi(\partial T\cap K)+{\rm M}_\varphi(S_2)\cr}$$
Thus
$$\eqalign{T(\omega)&\le {\rm M}_\varphi(\partial T\cap A)+T_1(\omega)+{\rm M}_\varphi(S_2)\cr
&\le {\rm M}_\varphi(\partial T\cap K)+{\rm M}_F(\partial T\setminus K).\cr}$$
There exists a continuous $\psi$ such that $\psi=\varphi$ on $A$, $\psi\ge\varphi$ on $K\setminus A$, and $\psi=F$ on $V\setminus K$. Then 
$T(\omega)\le {\rm M}_\psi(\partial T)$, thus Theorem 1 implies that there exists $\tau\in{\cal E}^q(V)$ with $d\tau=\omega$ and $||\tau||\le\psi+\varepsilon$. This implies Theorem 1'.

\vfill\eject
\noindent{\bf 3. The case of volume forms}
\medskip
{\bf Proof of theorem 2.}\hskip3mm Using the density of smooth currents and Theorem 1, it suffices to prove that
$T_h(\omega)\le {\rm M}_\varphi(\partial T_h)$
for every current of the form
$T_h(\varphi)=\int_V h\varphi$,
where $h$ is a smooth function with compact support. Then
$$V_\varphi(\partial T_h)=\sup\{\int_Vhd\tau\mid||\tau||\le f\}=\sup\{\int_Vdh\wedge\tau\mid||\tau||\le f\}=\int_V||dh||\varphi\ \nu$$
where $\nu$ is the Riemannian volume form.
\smallskip
By the coarea formula [F] applied to $|h|$,
$\int_V||dh||f\nu=\int_0^{+\infty}(\int_{|h|= t}\varphi d\sigma)\wedge dt.$
For almost all $t$, $\Omega_t=\{|h|\ge t\}$ is a smooth compact domain with boundary $\{{|h|= t}\}$. The hypothesis implies
$$V_\varphi(\partial T_h)\ge\int_0^{+\infty}(\int_{|h|\ge t}\omega)\wedge dt=\int_V|h|\omega.$$
Since this is $\ge T_h(\omega)$, Theorem 3 is proved.
\vfill\eject
\noindent{\bf 4. Filling and cofilling invariants}

\medskip
We recall here several definitions given by Gromov in [G2], chap.5 (and some variants).
\smallskip
Let $\gamma:S^1\to V$ be a rectifiable loop, homologous to zero. The {\it filling area} $\hbox{\rm Fill Area}(\gamma)$  is the infimum of the area of an integer singular $2$-chain $c$ with boundary $\gamma$. 
If we take the infimum over all real chains, we obtain the {\it real filling area}  $\R\hbox{\rm Fill Area}(\gamma)$, which is defined as soon as $\gamma$ is real-homologous to zero. If
$\gamma$ is integer-homologous to zero,
$\R\hbox{\rm Fill Area}(\gamma)=\lim_n\hbox{\rm Fill Area}(\gamma^n)/n.$
\smallskip
We can define analogously  $\R\hbox{\rm Fill Area}(b)$ for any real singular boundary. 
It clearly depends only on $I_b$. 
In fact, one can define (in any dimension) the {\it filling mass} of a boundary current:
$$\hbox{\rm Fill Mass}(S)=\inf\{{\rm M}(T)\mid T\in\partial{\cal E}'_q(V)\ {\rm and}\ \partial T=S\}.$$
By the density theorem, $\hbox{\rm Fill Mass}(I_b)=\R\hbox{\rm Fill Area}(b)$.
\medskip
The following result is proved in [F2], 4.13. Actually, it is only stated for locally flat forms, but regularization immediately gives the result with smooth forms (cf also [GLP], 4.35). 
\medskip

{\bf Whitney's duality.}\hskip3mm {\it If $S_0\in \partial{\cal E}'_q(V)$, 
$$\hbox{\rm Fill Mass}(S_0)
=\sup\{ S_0(\tau)\mid \tau\in {\cal E}^{q-1}(V)\ ,\ ||d\tau||\le1\}.$$}

\medskip
We recall the proof for the convenience of the reader. The argument is quite close to the proof of Theorem 1.
\medskip
The inequality $\ge$ is an immediate consequence of Stokes. To prove $\le$, we have to find for every $\varepsilon >0$ a smooth form $\tau$ such that $||d\tau||\le1$ and $S_0(\tau)\ge
{\rm Fill Mass}(S_0)-\varepsilon$. It suffices to find a locally flat form with these properties, then regularization will give the desired smooth one. 
\smallskip
Actually we can then take $\varepsilon=0$. Indeed, by Hahn-Banach there exists a linear form $t$ on ${\cal F}_{q-1}(V)$ such that
$t(S_0)=\hbox{\rm Fill Mass}(S_0)$ and $|t(S)|\le\hbox{\rm Fill Mass}(S)$ for every $S$ which is a boundary. This is equivalent to a flat form $\tau$ such that $S_0(\tau)=\hbox{\rm Fill Mass}(S_0)$ and $|S(\tau)|\le\hbox{\rm Fill Mass}(S)$ for every $S$ which is a boundary, which in turn is equivalent to: $|T(d\tau)|\le{\rm M}(T)$ for every $T$, ie $||d\tau||\le1$.

\medskip

{\bf Cofilling function.}\hskip3mm Fix $x_0$ in $V$. Gromov defines {\it the} cofilling function as `` the infimum of all functions'' $f:\R_+\to\R_+$ such that 
every exact $2$-form $\omega$ on $V$ with $||\omega||\le1$ has a primitive $\tau$ on $V$ satisfying $||\tau(x)||\le f(d(x_0,x))$. 

To make this more precise, we say that such a function $f$ is  {\it a} cofilling function, and we define
${\rm Cof}_q(R)$ as the infimum of all $C\ge0$ such that every exact $2$-form $\omega$ on $V$ with $||\omega||\le1$ has a primitive $\tau$ on $V$ satisfying $||\tau||\le C$ on $B'(x_0,R)$.
\smallskip
For $q=2$, we set ${\rm Cof}={\rm Cof}_2$.
Under reasonable assumptions, we shall see that $C{\rm Cof}(CR)$ is a cofilling function for some constant $C$, 
which will justifiy Gromov's  definition.
\medskip
We begin by a general geometric characterization of ${\rm Cof}_q$.

\medskip
{\bf Proposition 1.}
{\it For every $R\ge0$, 
$${\rm Cof}_q(R)=\sup\{{\hbox{\rm Fill Mass}(S)\over {\rm M}(S)}\mid S\in \partial{\cal E}'_q(V)\ , \ {\rm supp}(S)\subset  B(x_0,R)\}.$$}
\medskip
{\bf Proof.}\hskip3mm Using Whitney's duality, it suffices to prove that, for every $\omega\in d{\cal E}^q(V)$ with $||\omega||\le1$, one has
$$\inf_{\tau,d\tau=\omega} \max_{B'(x_0,R)}||\tau||=\sup\{{T(\omega)\over{\rm M}(\partial T)}\mid T\in{\cal E}'_q(V)\ , \ {\rm supp}(\partial T)\subset B(x_0,R)\}.$$
\smallskip
The inequality $\ge$ is obvious by Stokes. To prove $\le$, denote by ${\cal R}$ the right-hand-side.
We need to find, for every $\varepsilon>0$, a primitive $\tau$ with $||\tau||\le {\cal R}+\varepsilon$ on $B(x_0,R)$. 
This results from Theorem 1' with $U=B(x_0,R)$ and $\varphi\equiv {\cal R}$.

\medskip
Now we suppose $q=2$, and $H_1(V;\R)=0$.
\medskip
{\bf Proposition 2.}\hskip3mm
{\it If $H_1(V;\R)=0$,
$${\rm Cof}(R)=\sup\{{\R\hbox{\rm Fill Area}(\gamma)\over \ell(\gamma)}\mid \gamma\in {\rm Lip}(S^1,B(x_0,R))\}.$$}
{\it Remark.}\hskip3mm One may replace ${\rm Lip}$ by $C^\infty$.
\medskip
{\bf Proof.}\hskip3mm The hypothesis implies that $I_\gamma\in\partial{\cal E}'_2(V)$ for every Lipschitz loop $\gamma$. In Proposition 1, we may restrict by density and homogeneity to $S=\sum_{i=1}^k I_{\gamma_i}$ where $\gamma_i\in {\rm Lip}(S^1,B(x_0,R))$. 
\smallskip
If ${\rm M}(S)<\sum\ell(\gamma_i)$, the loops have common parts which cancel. By approximation,  we may assume that this common part is defined on unions of segments. By surgery, one has $S=\sum I_{\gamma'_j}$ with no cancellations.
Thus we may assume that ${\rm M}(S)=\sum\ell(\gamma_i)$
\smallskip
Then $\hbox{\rm Fill Mass}(S)=\R\hbox{\rm Fill Area}(\sum\gamma_i)\le\sum_{i=1}^k\R{\rm Fill Area} (\gamma_i)$ and thus
$${\hbox{\rm Fill Mass}(S)\over{\rm M}(S)}\le{\sum_{i=1}^k\R{\rm Fill Area}\over\sum\ell(\gamma_i)}\le{\max_i\R{\rm Fill Area}(\gamma_i)\over\ell(\gamma_i)}.$$
This proves Proposition 2.

\medskip
{\bf Real filling area function.}\hskip3mm This is the function $\R FA:\R_+\to\R_+$ defined by
$$\R FA(R)=\sup\R\hbox{\rm Fill Area}(T_c)\mid c\in {\rm Lip}(S^1,M)\ ,\ [c]=0\in H_1(V,\R)\ ,\ \ell(c)\le R\}.$$ 
Using $c^n$, one sees that $\R FA(nR)\ge n{\rm FA}(R)$ if $n\in\N$, thus $\displaystyle{\R FA(R)\over R}$ is ``almost non-decreasing'': $\displaystyle{\R FA(r)\over r}\le2\displaystyle{\R FA(R)\over R}$ if $r<R$.
 \medskip
{\bf Theorem 3.}\hskip3mm {\it (i) If  $H_1(V;\R)=0$,
$\displaystyle{\rm Cof}(R)\le 2{\R FA(3R)\over R}.$

(ii) If moreover $V$  is $C$-quasihomogeneous and $R\ge 2C$, 
 $\displaystyle{\rm Cof}(R)\ge C^{-3} {\R FA(R)\over R}.$
 Thus $\displaystyle{\rm Cof}(R)\sim {\R FA(R)\over R}$.}
\medskip
{\bf Proof}

\smallskip

(i)
Let $\gamma$ be a loop in $B(x_0,R)$. If $\ell(\gamma)\le R$,
$${\R{\rm Fill Area}(\gamma)\over\ell(\gamma)}\le\sup_{r\le R} {\R FA(r)\over r}\le 2{{\rm FA}(R)\over R}.$$

If $\ell(\gamma)> R$, we take $x_1,\cdots,x_k\in c$ with $k$  the smallest integer $\ge \ell(\gamma)/R$, such that the length of the arc $\gamma_i=x_ix_{i+1}$
on $\gamma$ is at most $R$, where we identify $x_{k+1}=x_1$. 
We define an oriented loop
$\gamma'_i=f_i*\gamma_i*f_{i+1}^{-1}$ where $f_i$ is a path from $x_0$ to $x_i$ of length $\le R$. 
Then $\ell(\gamma'_i)\le 3R$ and $I_\gamma=\sum_{i=1}^k I_{\gamma_i}$, thus
$$\R{\rm Fill Area}(\gamma)\le\sum_{i=1}^k\R FA(\gamma_i)\le k\R FA(3R)\le({\ell(\gamma)\over R}+1)\R FA(3R).$$
 Thus
$${\R{\rm Fill Area}(\gamma)\over\ell(\gamma)}\le{\R FA(3R)\over R}(1+{R\over \ell(\gamma)})\le2{\R FA(3R)\over R}.$$
Taking the supremum over all $\gamma$ and using Proposition 2, we obtain
(i).

\smallskip
(ii)
Let $\gamma\in {\rm Lip}(S^1,V)$ be a loop of length $\le R$. Its diameter is at most $R/2$, thus the quasihomogeneity gives $\gamma'=\varphi\circ\gamma$ with values in $B(x_0,R/2+C)$, of length $\le CR$. 
It also implies $\R{\rm Fill Area}(\gamma)\le C^2\R{\rm Fill Area}(\gamma')$.
For $R\ge 2C$, $\gamma'(S^1)\subset B(x_0,R)$, thus $\R{\rm Fill Area}(\gamma')\le {\rm Cof}(R)\ell(\gamma')\le CR {\rm Cof}(R)$. 
\smallskip Finally, 
$\R{\rm Fill Area}(\gamma)\le C^3R{\rm Cof}(R),$
which gives (ii).
\medskip
{\bf Proposition 3.}\hskip3mm {\it We make the assumptions (i) and (ii) of Theorem 3. Define $$\varphi(x)=3C^2\displaystyle{\R FA(6d(x_0,x))\over d(x_0,x)}.$$ Then
$\hbox{\rm Fill Mass}(S)\le {\rm M}_{\varphi}(S)$
for every $S\in \partial{\cal E}'_2(V)$.}
\medskip
{\bf Proof.}\hskip3mm As in Proposition 2, we first  reduce to the case where $S=I_\gamma$ with $\gamma$ a loop. Then using the quasihomogeneity, we may assume 
that $$\gamma(S^1)\subset B(x_0,\ell(\gamma)/2+C)\setminus B'(x_0,1)\subset B(x_0,\ell(\gamma))\setminus B'(x_0,1),$$
and also that $\gamma(0)\in B(x_0,C)$. This will increase the constant $N$ by at most a factor $C^2$.
\smallskip
We may assume that $\gamma:[0,\ell(\gamma)]\to V$ is parametrized by arclength. 
We define
$t_0=0$ and $t_i=t_{i-1}+{1\over 2}d(x_0,\gamma(t_{i-1}))$ as long as $t_i\le\ell(\gamma)$. Since $d(x_0,\gamma(t_i))\ge1$, this is possible up to a maximal $i=N$. We obtain thus  $N$ consecutive arcs $I_i=\gamma|[t_{i-1},t_i]$, $1\le i\le N$.
\smallskip
Let $c_i$ be a minimal  geodesic from $x_0$ to $\gamma(t_i)$, and let $\gamma_i$ the loop $c_{i-1}*(\gamma|I_i)*c_i^{-1}$. Define also $\gamma_0=c_N*\gamma|[t_N,\ell(\gamma)]*c_0^{-1}$.
Set 
$$\eqalign{d_i&=d(x_0,\gamma(t_i))\ , \ \ell_i=\ell(\gamma_i)=t_i-t_{i-1}\ , \ \ell_0=\ell(\gamma)-t_N\cr
\delta_i&=d(x_0,\gamma(I_i))\ ,\ \Delta_i=\max_{t\in I_i}d(x_0,\gamma(t)).\cr}$$ 
Then $\ell_i={1\over2}d_{i-1}$, thus
$$\eqalign{\delta_i&\ge d_{i-1}-\ell_i=\ell_i\cr
d_i\le \Delta_i&\le d_{i-1}+\ell_i= 3\ell_i\le 3\delta_i.\cr}$$
Thus
$$\eqalign{\R{\rm Fill Area}(\gamma)&\le\sum_{i=0}^N\R{\rm Fill Area}(\gamma_i)\le\R FA(d_N+\ell_0+d_0)+\sum_{i=1}^N\R FAd_{i-1}+\ell_i+d_i)\cr
&\le\R FA(d_N+\ell_0+d_0)+\sum_{i=1}^N{3\R FA(6\delta_i)\ell_i\over\Delta_i}\cr}$$
Moreover, $d_0\le C$, $\ell_0<{1\over2}d_N$ by maximality. Thus 
$$d_N\le d(\gamma(t_0),\gamma(t_N)) +C\le \ell_0+C<{1\over2}d_N+C,$$
so $d_N<2C$, $d_N+\ell_0+d_0\le 2C+C+C=4C$. 

Defining $\displaystyle\psi(x)={3\R FA(6d(x_0,x))\over d(x_0,x)}$, we have
$$\eqalign{\hbox{\rm Fill Mass}(I_\gamma)=\R{\rm Fill Area}(\gamma)&\le\R FA(4C)+\sum_{i=1}^N\min_{t\in [t_{i-1},t_i]}\psi(\gamma(t))\ (t_i-t_{i-1})\cr
&\le\R FA(4C)+\sum_{i=1}^N\int_{t_{i-1}}^{t_i}\psi(\gamma(t))\ dt\cr
&=\R FA(4C)+\int_0^\ell\psi(\gamma(t))\ dt\cr
&=\R FA(4C)+{\rm M}_\psi(I_\gamma).\cr}$$
Replacing $\gamma$ by $\gamma^n$ and making $n\to+\infty$, we deduce
${\rm Fill Mass}(I_\gamma)\le {\rm M}_\psi(I_\gamma).$
Thus for every $S\in\partial{\cal E}_2'(V)$, we have
$\hbox{\rm Fill Mass}(S)\le C^2{\rm M}_\psi(I_\gamma).$ This proves Proposition 3.

\medskip
{\bf Corollary.}\hskip3mm {\it Assume that $H_1(V;\R)=0$ and that $V$ is $C$-quasihomogeneous. Then every exact $2$-form with norm $\le1$ has a primitive such that $||\tau(x)||\le 4C^2\displaystyle{\R FA(6d(x_0,x))\over d(x_0,x)}.$ In other words, $f(x)=4C^2\displaystyle{\R FA(6d(x_0,x))\over d(x_0,x)}$ is a cofilling function.}
\medskip
By Theorem 3,(ii), it is the ``smallest'' cofilling function up to equivalence.

\vfill\eject
\noindent{\bf 5. Primitives of cocycles  of degree $2$ on a group}
\medskip
Recall that a $q$-cochain $u\in C^q(G;\R)$ on the group $G$ is a function $u:G^{q+1}\to\R$. The differential is defined by
$$du(g_0,\cdots,g_{q+1})=\sum_{i=0}^{q+1}(-1)^iu(g_0,\cdots,\widehat g_i,\cdots,g_q).$$
Recall that the subcomplex of $G$-invariants cochains  $C^*_{inv}(G;\R)$ gives rise to the group cohomology $H^*(G,\R)$.

\medskip
Recall
the statement of Theorem 4.
\medskip
{\bf Theorem 4.}\hskip3mm {\it Let $b$ be a $2$-cocycle on $G$, and let $F$ be a function from $G$ to $\R_+$. 
Then the following are equivalent:
\smallskip
(i) There exists $t\in C^1(G;\R)$ 
such that $da=b$ and
$ |a(g,g\overline s)|\le F(g).$
\smallskip
(ii) For every $g\in G$ and every relation $w=s_1^{\varepsilon_1}\cdots s_n^{\varepsilon_n}\in R$, one has, setting $g_i=g\overline s_1^{\varepsilon_1}\cdots\overline s_{i-1}^{\varepsilon_{i-1}}$ ($g_0=1$):
$$\big|\displaystyle\sum_{i=1}^nb(1,g_i,g_{i+1})-\sum_{\varepsilon_i=-1} b(1,g_{i+1},g_i)\big|
\le\sum_{i=1}^nF(g_i).$$
}
Using the canonical primitive $a_0(g,h)=b(g_0,g_1)$, we can write $a=a_0+dm$ with $m:G\to\R$ (ie $a(g_0,g_1)=a_0(g_0,g_1)+m(g_1)-m(g_0)$). Setting then $\alpha_0(g,s)=a_0(g,g\overline s)$, we see that the significant data is $\alpha_0:G\times S\to\R$, which we can view as function on the edges of the Cayley graph. We can restate Theorem 4 as follows.
\medskip
{\bf Theorem 4'.}\hskip3mm {\it Let $\alpha_0$ be a function on $G\times S$, and let $F$ be a function from $G$ to $\R_+$. 
Then the following are equivalent:
\smallskip
(i) There exists $m:G\to\R$ 
such that $ |\alpha_0(g,s)+m(g\overline s)-m(g)|\le F(g)$.
\smallskip
(ii) For every $g\in G$ and every relation $w=s_1^{\varepsilon_1}\cdots s_n^{\varepsilon_n}\in R$, one has, setting $g_i=g\overline s_1^{\varepsilon_1}\cdots\overline s_{i-1}^{\varepsilon_{i-1}}$ ($g_0=1$):
$$\big|\displaystyle\sum_{\varepsilon_i=1}\alpha_0(g_i,s_i)-\sum_{\varepsilon_i=-1} \alpha_0(g_{i+1},s_i)|
\le\sum_{i=1}^nF(g_i)$$
}
\medskip
{\bf Proof of Theorem 4'.}\hskip3mm We consider $m$ as a linear form on $\R[G]$. By Hahn-Banach, (i) is equivalent to
$$\sum_{i=1}^n\tau_i(g_i\overline s_i-g_i)=0 \Rightarrow\big|\sum_{i=1}^n\tau_i\alpha_0(g_i,s_i)\big|\le\sum_{i=1}^n|\tau_i|F(g_i),\leqno(i)'$$
where the $\tau_i$ are nonzero real numbers.
\medskip
1) Suppose that (i) is true. The hypothesis of (ii) implies 
$\sum_{i=1}^n\ \varepsilon_i(g_i\overline s_i^{\varepsilon_i}-g_i)=0.$
Then (i)' with $\tau_i=\varepsilon_i$ gives (ii).
\medskip
2) Suppose that (ii) is true, and that 
$$\sum_{i=1}^n\tau_i(g_i\overline s_i-g_i)=0.\leqno(1)
$$
We want to prove that $$\big|\sum_{i=1}^n\tau_i\alpha_0(g_i,s_i)\big|\le\sum_{i=1}^n|\tau_i|F(g_i).\leqno(2)$$
\smallskip
We argue by induction over $n$, the result being trivial for $n=0$. We may assume that $\tau_1>0$ and that $|\tau_1|$ is minimal.
\smallskip
The term $\tau_1g_1\overline s_1$ must cancel with some other, ie there exists $i=i_2$ such that either ($g_1\overline s_1=g_i$ with $\tau_i\tau_1>0$), or ($g_1\overline s_1=g_i\overline s_i$ with $\tau_i\tau_1<0$). 
Continuing with the term $\tau_ig_i\overline s_i$ or $\tau_ig_i$ respectively, we define inductively $i_1=1,i_2,i_3,\cdots$ and $\varepsilon_1=1,\varepsilon_2,\cdots$, such that, for all $k$, one has
$$\eqalign{g_1\overline s_{i_1}^{\varepsilon_1}\overline s_{i_2}^{\varepsilon_2}\overline s_{i_3}^{\varepsilon_3}\cdots\overline  s_{i_k}^{\varepsilon_k}&=\left\{\eqalign{g_{i_{k+1}}\ &\  {\rm if}\ \varepsilon_{k+1}=1\cr
g_{i_{k+1}}\overline s_{i_{k+1}}\ &\  {\rm if}\ \varepsilon_{k+1}=-1\cr}\right.\cr
\varepsilon_k&={\rm sgn}(\tau_{i_k}).\cr}$$
Let $k$ be the smallest integer such that $i_{k+1}=i_1=1$. If we have $i_\ell=i_m$ for some $1\le\ell<m\le k$, we can suppress the indexes $i_r$ with $r$ between $\ell+1$ and $m$. Thus we can assume that all the $i_r$ are distinct. 
Since $g_{k+1}=g_1$, we have $s_{i_1}^{\varepsilon_1}\cdots s_{i_k}^{\varepsilon_k}\in R$, with $\epsilon_1=1$. This implies
$$\sum_{\varepsilon_r=1}(g_{i_r}\overline s_{i_r}-g_{i_r})-
\sum_{\varepsilon_r=-1}(g_{i_{r+1}}\overline s_{i_r}-g_{i_{r+1}})=0.$$
Changing the numbering of the $g_i$, we can rewrite this equality and $(1)$ as
$$\eqalign{\sum_{i=1}^k\varepsilon_i(g_i\overline s_i-g_i)&=0\cr
\sum_{i=1}^n\tau_i(g_i\overline s_i-g_i)&=0.\cr}$$
We also have $\varepsilon_i={\rm sgn}(\tau_i)$.
Combining the two, we get
$$\sum_{i=2}^k(\tau_i-\varepsilon_i\tau_1)(g_i\overline s_i-g_i)+\sum_{i=k+1}^n=0.$$
The inductive hypothesis implies
$$\Big|\sum_{i=2}^k(\tau_i-\varepsilon_i\tau_1)\alpha_0(g_i,s_i)+\sum_{i=k+1}^n\tau_i\alpha_0(g_i,s_i)\Big|\le \sum_{i=2}^k|\tau_i-\varepsilon_i\tau_1|F(g_i)+\sum_{i=k+1}^n|\tau_i|F(g_i).$$
By (ii), the property $s_{i_1}^{\varepsilon_1}\cdots s_{i_k}^{\varepsilon_k}\in R$, with $\epsilon_1=1$ implies (with the new numbering)
$$\big|\sum_{\varepsilon_i=1}\varepsilon_i\alpha_0(g_i,s_i)\big|\le\sum_{r=1}^kF(g_{i_r}).$$
Finally, the hypotheses imply $|\tau_i-\varepsilon\tau_1|+|\tau_1|=|\tau_i|$, thus combining the last two inequalities gives  (i)'. This finishes the proof of Theorem 4.
\medskip
{\bf Remark.}\hskip3mm The proof of $(ii)\Rightarrow(i)'$ is related to the property $\ker(\partial_1)={\rm im}(\theta)$ in the ``Hopf'' exact sequence
$$0\rightarrow R^{ab}\hfl[-1mm][8mm][-1mm]{\theta}{}\Z[G]^p
\hfl[-1mm][8mm][-1mm]{\partial_1}{}\Z[G]\hfl[-1mm][8mm][-1mm]{\varepsilon}{}\Z\rightarrow0,$$
or its tensorization over the reals. Here the relation module $R_{ab}$ is the abelianization of $R\subset F(S)=F_p$, the relation subgroup. The $G$-action comes from conjugation in $F_p$, thus $\theta([gwg^{-1}])=g\theta([w])$.
\smallskip
Define $$x=\sum_{i=1}^n\tau_ig_ie_{s_i}\in\R[G]^S=\R[G]^p.$$
The relation  
$\sum_{i=1}^n\tau_i(g_i\overline s_i-g_i)=0$ translates to $\partial_1(x)=0$. Thus there exists a finite family $(r_q,\mu_q)\in R\times \R$ such that $$\theta(\sum_q\mu_q\overline r_q)=\sum\tau_ig_i.\leqno(4)$$
By [Brow] p.45, an explicit formula for $\theta([w])$, where $r$ is the relation $s_1^{\varepsilon_1}\cdots s_k^{\varepsilon_k}$,  is
$$\theta([w])=\sum_{s\in S}\overline{\partial r\over\partial s}\ e_s=\sum_{\varepsilon_i=1}g_ie_{s_i}-\sum_{\varepsilon_i=-1}g_{i+1}e_{s_i},$$
where $g_i=\overline s_1^{\varepsilon_1}\cdots\overline s_{i-1}^{\varepsilon_{i-1}}$.
\smallskip
In other words, $\theta$ is induced by the derivation $d:F\to\Z[G]^p$ such that $d(s_i)=e_i$.
\smallskip

Then the formula $(*)$ translates into a decomposition of the identity $\sum\tau_i(g_i\overline s_i-g_i)=0$ into a combination of identities ($\sum \varepsilon_i(\widetilde g_i-\widetilde g_i\overline s_i)=0$) associated to the relations.

\vfill\eject
\noindent{\bf 6. Relation with the $\ell_1$-norm of Gersten and the homological Dehn function}

\medskip
Assume now that  $G=\langle s_1,\cdots,s_p\mid r_1,\cdots,r_q\rangle$ be a finitely presented group.
Consider the exact sequence associated to the cellular homology of $\widetilde M$, where $M$ is the $2$-complex defined by the presentation: 
$$\Z[G]^q\hfl[-1mm][8mm][-1mm]{\partial_2}{}\Z[G]^p
\hfl[-1mm][8mm][-1mm]{\partial_1}{}\Z[G]\hfl[-1mm][8mm][-1mm]{\varepsilon}{}\Z\rightarrow0,$$
or its tensorization over the reals. The Hopf exact sequence gives an isomorphism $\theta:R_{ab}\simeq\ker\partial_1\subset\Z[G]^p$ (see the previous section), and we have
$\theta([r_i])=\partial_2(f_i)$ for $1\le j\le q$.
\medskip
The group $\Z[G]$ (or the vector space $\R[G]$) is equipped with the $\ell_1$-norm
$|\sum_g\tau_g g|_1=\sum_g|\tau_g|$. This extends to $\Z[G]^q$, $\Z[G]^p$.
Then one can define another norm on $\ker\partial_1={\rm im}(\partial_2)$:
$$||z||=\inf\{|c|_1\mid c\in\Z[G]^q\ , \ \partial_2c=z\}$$
(If we work with integer coefficients, we have a minimum).
\smallskip
If
 $w\in R$ is a relation, $[w]\in R^{ab}={\rm im}(\partial_2)$. S. Gersten in [Gersten 1990] gives the following definition: 
$$||[w]||=\inf\{|c|_1\mid c\in \Z[G]^q\ , \ \partial_2c=[w]\}.$$
One checks that, if the coefficients are integers, this is equal to the abelianized isoperimetric function of [BMS]:
$$||[w]||=\Delta^{ab}(w)=\min\{m\mid w\in\prod_{i=1}^m u_ir_{j_i}^{\varepsilon_i}u_i^{-1}[R,R]\}.$$

We shall need the stable version
$$\Delta^{ab}_{\hbox{\bbsmall R}}(w)=\lim_{n\to\infty}{\Delta^{ab}(w^n)\over n}.$$

\smallskip
{\bf $1$-cycle associated to a relation.}\hskip3mm Here it suffices that $G=F(S)/R$ be finitely generated. The space of $k$-chains is 
$C_k(G)=\Z[G^{k+1}]$. As a $\Z[G]$-module, it free with the standard basis
$$[g_1|\cdots|g_k]=(1,g_1,g_1g_2,\cdots,g_1g_2\cdots,g_n).$$
\smallskip

If $w=s_1^{\varepsilon_1}\cdots s_n^{\varepsilon_n}\in R$, we define $g_i=\overline s_1^{\varepsilon_1}\cdots\overline s_{i-1}^{\varepsilon_{i-1}}$ and
$$I_w=\sum_{i=1}^n(g_i,g_{i+1})=\sum_{i=1}^n g_i[\overline s_i]\in C_1(G).$$
In other words, $I_w=\eta(\theta([w]))$, where $\eta:\Z[G]^p\to C_1(G)$ is $\R[G]$-linear and 
$\eta(e_i)=[\overline s_i]$.
\smallskip
Clearly, $I_w$ is a cycle, ie $I_w\in Z_1(G;\R)$, and $I_w$ only depends on $[w]\in R^{ab}$.
Then one has simply
$I_w=\sum_{i=1}^n[g_i,g_{i+1}].$
\smallskip
The complex $C_*(G;\R)$ is exact, thus there exists $T\in C_2(G;\R)$ with $\partial T=I_w$.
\smallskip
{\bf Proposition.}\hskip3mm {\it The map $[w]\mapsto I_w$ is injective from $R^{ab}$ to $Z_1(G)$.}
\medskip
{\bf Proof.}\hskip3mm  View $w$ as a loop starting from $1$ in the Cayley graph
of $(G,S)$. The property $I_w=0$ means that $w$ has an algebraic coefficient $1$ on each edge. This means that it is homologous to zero, ie $w\in[R,R]$, or $[w]=0$, qed.
\medskip
{\bf Question.}\hskip3mm O\`u y a-t-il une r\'ef\'erence \`a \c ca dans la litt\'erature ?
\medskip
{\bf Corollary.}\hskip3mm {\it If  $w\in R$,
$$\Delta^{ab}_{\hbox{\bbsmall R}}(w)=\max\{a(I_w)\mid t\in C^1(G;\R)\ {\rm and}\ (\forall (g, j)\ |a(gI_{r_j})|\le1\}.$$
}
\smallskip
{\bf Remark.}\hskip3mmNote the similarity with (i) in the lemma of section 4.
\medskip
{\bf Proof of the corollary.}\hskip3mm Let $w\in R$. By the lemma,
$$w\in\prod_{i=1}^m u_ir_{j_i}^{\varepsilon_i}u_i^{-1}[R,R]\Leftrightarrow I_w=\sum _{i=1}^m\varepsilon_ig_iI_{r_{j_i}}.$$
Thus $\Delta^{ab}(w)=\min\{m\in\N\mid I_w=\sum_{i=1}^m\varepsilon_ig_iI_{r_{j_i}}\}$, which implies
$$\Delta^{ab}_{\hbox{\bbsmall R}}(w)=\inf\{\sum|\tau_i|\mid I_w=\sum \tau_ig_iI_{r_{j_i}}\},$$
the sums being finite and with real coefficients.
The corollary is an immediate consequence of Hahn-Banach.
\vfill\eject
\noindent{\bf 7. Filling and cofilling in groups}

\medskip
Here $G=\langle s_1,\cdots,s_p\mid r_1,\cdots,r_q\rangle$ is a group equipped with a finite presentation. 
This gives a norm function for each $2$-cocycle $b\in Z_2(G)$: if $b=da$, one sets
$$||b||(g)=\max_j|a(gI_{r_j})|.$$
Since $I_w$ is closed, it is a boundary $I_w=\partial_2(c_w)$, thus $a(gI_{r_j})=b(gc_{r_j})$ depends only on $b$.
\medskip
{\bf Cofilling function.}\hskip3mm For $n\in\N$, we define
${\rm Cof}(n)$ as the infimum of all $C\ge0$ such that every cocycle $b$ on $G$ with $||b||\le1$ has a primitive $a$ satisfying $||u_a||\le C$
on $B_S(n)$, ie $|a(g,g\overline s^\pm1)|\le C$ if $|g|\le n$. 
\medskip
{\bf Lemma.}\hskip3mm {\it For every $n\in\N$, one has
$${\rm Cof}(n)=\sup\{{\Delta^{ab}_{\hbox{\bbsmall R}}(w)\over |w|}\mid w\in R\ ,\ |w|\le n\}.$$
}
{\bf Proof.}\hskip3mm Recall the corollary in section 6:
$$\Delta^{ab}_{\hbox{\bbsmall R}}(w)=\max\{a(I_w)\mid t\in C^1(G;\R)\ , \ ||da||\le1\}.$$ 
Thus it suffices to prove that, for every $b\in dC^1(G;\R)$ with $||b||\le1$, one has
$$\inf_{a,da=b} \max_{B_S(n)}||a||=\sup\{{b(T)\over|\partial T|_1}\mid T\in C_2(B_S(n))\}.$$

\smallskip
Call ${\cal L}$ the left-hand side and ${\cal R}$ the right-hand side. The inequality $({\cal R}\le {\cal L})$ is obvious by Stokes.
To prove that $({\cal L}\le {\cal R})$,  we need to find  a primitive $a$ with $||a||\le {\cal R}$ on $B_S(n)$. 
For this, we apply Theorem 4 with $F={\cal R}$ on $B_S(n)$ and $F=\infty$ elsewhere.
\smallskip
It suffices to have
$|b(T)|\le {\cal R}|\partial T|_1$ for every $T\in C_2(B_S(n))$.
We have $\partial T=\sum\tau_iI_{w_i}$ with $|\partial T|_1=\sum|\tau_i||w_i|$, thus we may assume 
$\partial T=I_w$. Then
$$|b(T)|\le \Delta^{ab}_{\hbox{\bbsmall R}}(w)\le {\cal R} |I_w|_1={\cal R}|\partial T|_1,$$
which proves the lemma.
\medskip
Thus we obtain  the {\it homological Dehn function}, or {\it abelian isoperimetric function} [BMS]:
$$\delta^{ab}(n)=\sup\{\Delta^{ab}(w)\mid |w|\le n\}=\sup\{||z||\mid z\in\ker\partial_1\ , \ |z|_1\le n\}.$$
Again, there are two versions, with integer or real coefficients.
\smallskip
Let $w=s_1^{\varepsilon_1}\cdots s_n^{\varepsilon_n}$ be a relation, and let
$w=\prod_{k=1}^Nx_kr_{j_k}^{\eta_k}x_k^{-1}$
be a decomposition into elementary relations modulo $[R,R]$. Assume that $N$ is minimum, ie $N=\Delta^{ab}(w)$.
Then
$\theta([w])=\sum_{k=1}^N\eta_k\theta(x_k[r_{j_k}]),$ and
the condition on $b$ to have a primitive bounded by $F$ becomes 
$$|\sum_{k=1}^N\eta_kb(x_kc_{r_{j_k}})|\le\sum_{\varepsilon_i=1}F(g_i)+\sum_{\varepsilon_i=-1}F(g_{i+1}).$$
Let $M=\max(|[c_{r_j}]|)$, then the left-hand-side is bounded by $M\Delta^{ab}(w)$.
Replacing $w$ by $w^n$ and making $n\to+\infty$, we see that it is in fact bounded by $\Delta^{ab}_{\bbsmall R}(w)$.

\medskip
{\bf Primitive of a bounded cocycle.}\hskip3mm 
In order for Question 2 to have a positive answer for every $a\in Z^2G$ with $||a|||\le1$, it suffices that, for every relation $w\in R$, one have
$$\Delta^{ab}_{\bbsmall R}(w)\le M^{-1}\big(\sum_{i=1}^nF(g_i)\big).$$

\bigskip
{\bf Special case: constant bounds.}\hskip3mm Let 
$f=A$ be constant in Question 2. 
Then Theorem 4 says that the answer is positive if, for every relation $w\in R$, one has
$\Delta^{ab}_{\bbsmall R}(w)\le AM^{-1}|w|,$ ie $\delta^{ab}_{\bbsmall R}(n)\le AM^{-1}n$.
\medskip
{\bf Relation with hyperbolicity.}\hskip3mm
By Mineyev, this is equivalent to the hyperbolicity of $G$.
\bigskip
{\bf Primitives of cocycles of degree $>2$ on a group}

\medskip
Let $q$ be an integer $>2$. Let $G$ be a group of type $F_q$, ie there exists a finite cell complex $M$ such that $\pi_1M=G$ and $\pi_iM=0$ for 
$2\le i\le q-1$. Alternatively, there exists a cell complex $Y$ which is a $K(G,1)$ and has a finite $q$-skeleton. 
\vfill\eject
\noindent{\bf 8. Relation between Questions 1 and 2}
\bigskip
Let $V$ be a Riemannian manifold equipped with a geometrically bounded triangulation $T$. Let $I^*:C^*(T)\to{\cal E}^*(V)$ be the integration morphism. The following result is contained in substance in [Si]. The proof consists  essentially in adding bounds to the proof of the theorem of de Rham given in [ST], p.165 sqq.\medskip
{\bf Proposition}\smallskip {\it 
(i) There exists  a right inverse $R^*$ for $I^*$ which is a chain map ($R\circ d=d\circ R$) and satisfies
$$||R(u)_x||+||d(R(u))_x||\le C\max\{|u(\sigma)|\mid \sigma\subset B'(x,C)\}.$$
\smallskip
(ii)There exists a linear map $\Pi^q:{\cal B}^q(V)\cap\ker(I^q)\to{\cal E}^{q-1}(V)$ such that $\Pi^q(\omega)$ is a primitive of $\omega$ and 
$$||\Pi^q(\omega)_x||\le C\max\{||\omega_y||\mid y\in B'(x,C)\}.$$}
 
 Actually the right inverse has been defined by Whitney ([W] p.226), the new observation is (ii). In fact, a stronger and more natural property holds.
 \medskip
 {\bf Proposition.}\hskip3mm {\it There exists a chain homotopy $H^*:R^*I^*-{\rm Id}\simeq 0$, ie a linear map $H^*=(H^q:{\cal E}^q(V)\to {\cal E}^{q+1}(V))$ of degree $1$, with the property
 $$||H(\omega)_x||\le C\max\{||\omega_y||,||d\omega_y||\mid y\in  B'(x,C)\}.$$}

 \medskip
{\bf Corollary.}\hskip3mm {\it Let $\omega\in{\cal E}^q$ be an exact $q$-form on $V$, and let $t\in C^{q-1}(T;\R)$ be a primitive of $I^q(\omega)$. Then $\omega$ has a primitive $\tau\in{\cal E}^{q-1}(V)$ such that
$$||\tau_x||\le C(\max_{B'(x,C)}||\omega||+\max\{|t(\sigma)|\mid \sigma\subset B'(x,C)\}).$$}

{\bf Proof of the corollary.}\hskip3mm  Let $\omega_1=\omega-dR(t)$, so that $I^q\omega_1)=0$, and $\tau=d(R(t))+\Pi(\omega_1)$. Then $d\tau=\omega$, and the  
estimates are immediate.
\vfill\eject
\noindent{\bf 9. Relation between Questions 1, 2 and 3 for $q=2$}
\medskip
Let $M$ be a compact Riemannian manifold with infinite fundamental group, and 
$\pi:\widetilde M\to M$ be its
universal covering. 
\smallskip
Let $T$ be a smooth triangulation of $M$, which we lift to $\widetilde M$. We associate to $\omega$ the $2$-cochain $I_T(\omega)$.
\smallskip
Let $X$ be a smooth cellulation of $M$, with only one $0$-cell $x_0$. Thus $X^{(2)}$ defines a presentation of $\pi_1(M,x_0)=G$. Similarly, we lift $X$  to $\widetilde M$ and define the $2$-cochain $I_X(\omega$.
\smallskip
We have an action an action of $G=\pi_1(M,x_0)$ on $\widetilde M$. 
For each $g\in G$ choose a cellular path $\sigma(g)$ from $\widetilde x_0$ to $g\widetilde x_0$ representing $g$. This is the same as a normal form $\nu:G\to F$.
\smallskip

Let $\omega$ be an exact $2$-form on $\widetilde M$ for some $q\ge2$.
We define a $2$-cocycle $u\in C^2(G,\R)$ by
$$u(g_0,g_1,g_2)=\int_{D(g_0,g_1,g_2)}\omega,$$
where $D(g_0,g_1,g_2)\subset\widetilde M$ is any cellular disk [$C^1$ map defined on $D^2$] bounded by the loop
$$\gamma(g_0,g_1,g_2)=g_0(\sigma(g_0^{-1}g_1)*\sigma(g_1^{-1}g_2)*\sigma(g_0^{-1}g_2)^{-1}.$$
This is well defined since $\int_\Sigma\omega=0$ for every $2$-sphere $\Sigma\subset\widetilde M$. [in fact for any surface]

\smallskip
A primitive of $u$ is
$t_0(g_0,g_1)=u(1,g_0,g_1)=\int_{D(g_0,g_1)}\omega$
where $D(g_0,g_1)$ is any disk bounded by the loop
$\gamma(g_0,g_1)=\sigma(g_0)*\sigma(g_0^{-1}g_1)*\sigma(g_1)^{-1}.$

\medskip
We want to relate the following properties:
\smallskip
(1) There exists $\tau\in{\cal E}^1(\widetilde M)$ such that $d\tau=\omega$ and $||\tau||\le\varphi$.
\smallskip
(2) There exists $t\in C^1(\widetilde T)$ such that $dt=I_X(\omega)$ and $|t|\le f$.
\smallskip
(2') There exists $t\in C^1(\widetilde X)$ such that $dt=I_T(\omega)$ and $|t|\le f$.
\smallskip
(3) There exists $a\in C^1(G)$ such that $da=b$ and $|t(g,g\overline s^{\pm1})|\le F(g)$.
\medskip
{\bf Proposition}\hskip3mm
\smallskip
 {\it  (i) If (1) holds for some $\varphi$, (2) holds for
$f(\sigma)=C\max\{\varphi(x)\mid \sigma\subset B'(x,C)\}.$

\smallskip
(ii) If (2) holds for $f$, (1) holds for 
$\varphi(x)=C\max\{f(\sigma)\mid \sigma\subset B'(x,C)\}.$
\smallskip
(iii) If (2) holds for $f$, (3) holds for
$F(g)=C\max\{f(\sigma)\mid \sigma\subset {\rm st}^2(g\widetilde x_0)\}.$
\smallskip
(iii) If (3) holds for $F$, (2) holds for
$f(\sigma)=C\max\{F(g)\mid \sigma\subset {\rm st}^2(g\widetilde x_0)\}.$
}
\medskip
{\bf Proof.}\hskip3mm (i) is obvious: it suffices to take $t=I^1(\tau)$.

\smallskip
(ii) is an immediate consequence of the corollary in section 8.
\smallskip
(iii) and (iv). One defines $G$-equivariant chain maps $\psi_*:C_*(\widetilde X)\to C_*(G)$ and $\chi_*:C_*(G)\to C_*(\widetilde X)$ in degrees $\le2$ (cf. [Brown], p.46):
\smallskip
1) If instead of a triangulation we have a cellulation with $c_0=1$, $c_1=p$, $c_2=q$, we define 
\smallskip
\hskip5mm - $\psi_0=\chi_0={\rm Id}$;
\smallskip
\hskip5mm -  $\psi_1=\eta$, ie $\psi(e_i)=[\overline s_i]$, $1\le i\le p$
(cf section 6);
\smallskip
\hskip5mm
- for each $g\in G$, choose a normal form $\nu(g)=s_1^{\varepsilon_1}\cdots s_n^{\varepsilon_n}$ representing $g$, and set   

$$\eqalign{\chi_1([1,g])&=d(\nu(g))=\sum_{i=1}^n{\partial(\nu(g))\over\partial s_i}\cr
&=\sum_{\varepsilon_i=1}\ g_i e_{s_i}-\sum_{\varepsilon_i=-1}g_i e_{s_i},\cr}$$
where $g_i=\overline s_1^{\varepsilon_1}\cdots\overline s_{i-1}^{\varepsilon_{i-1}}$ as usual.
\smallskip
\hskip5mm - $\psi_2(f_j)={\rm Cone}(I_{r_j})$, $1\le j\le q$,  where ${\rm Cone}(g,h)=(1,g,h)=[g|g^{-1}h]$; if $g_1,\cdots,g_n$ are associated to $r_j$ as usual, $\psi_2(f_j)=\sum_{i=1}^n[g_i|\overline s_i^{\varepsilon_i}]$.
\smallskip
\hskip5mm - for each 
$(g,h)\in G\times H$, we choose a decomposition
$$\nu(g)\nu(h)\nu(gh)^{-1}=\prod_kx_kr_{j_k}^{\varepsilon_k}x_k^{-1},$$
or 
$$\nu(g)\nu(h)\nu(gh)^{-1}\equiv\prod_kx_kr_{j_k}^{\varepsilon_k}x_k^{-1}\ {\rm mod}\ [R,R].$$
Then we set
$$\chi_2([g|h])=\sum \varepsilon_k\overline x_k\sigma_{j_k}.$$
\medskip

\medskip
By duality we have cochain maps $\psi^*$ and $\chi^*$. Then

\bigskip
{\bf Relation between the three questions for $q>2$}
\medskip
We assume that $\pi^*\omega\in H^q(\widetilde M;\R)$ vanishes, ie there exists $u\in H^q(\pi_1M;\R)$ (unique) such that $i^*[u]=[\omega]$ where
$i:M\to X(\pi_1M,1)$ is the natural map (defined up to homotopy). 

\medskip
Assume that $\pi_1V$ is of type $F_q$ 
[or $\pi_iV=0$ for $2\le i\le q-1$]. 
\medskip

\vfill\eject
\centerline{\bf Bibliography}
\medskip
\item{[AG]} D. Allcock, S. Gersten, {\it A homological characterization of hyperbolic groups}, Invent. Math. 135 (1999), 723-742.
\smallskip
\item{[BMS]} G. Baumslag, C.F. Miller III, and H. Short, {\it Isoperimetric inequalities and the homology of groups}, 
Invent. Math. 113 (1993), 531-560.
\smallskip
\item{[Broo]} R. Brooks, {\it Some Riemannian and dynamical invariants of foliations}, in {\it Differential Geometry},
Birkh\"auser PIM 32, 1983, 56-72.
\item{[Brow]} K.S. Brown, {\it Cohomology of groups}, Springer, 1982.
\smallskip
\item{[F1]} W. Federer, {\it Geometric measure theory}, Springer GMW 153, 1969.
\smallskip
\item{[F2]} W. Federer, {\it Real flat chains, cochains and variational problems}, Indiana Math J. 24 (1974), 351-407.
\smallskip
\item{[Ge1]} S. Gersten, {\it Dehn functions and $\ell_1$-norms of finite presentations}, Algorithms and classification in combinatorial group theory, (G. Baumslag and C.F. Millett III ed), MSRI 8, Springer 1987, pp. 195-224.
\smallskip
\item{[Ge2]} S. Gersten, {\it A cohomological characterization of hyperbolic groups}, preprint 1996, available at http://math.utah.edu/$\sim$gersten.
\smallskip
\item{[Ge3]} S. Gersten, {\it Homological Dehn functions and the word problem}, preprint 1999, available as [G2].
\smallskip
\item{[Gr1]} M. Gromov, {\it Hyperbolic manifolds, groups and actions}, in {\it Riemannian surfaces and related topics (Stony Brooks 1979)},
Princeton Ann Math Studies 97, 1981.
\smallskip
\item{[Gr2]} M. Gromov, {\it K\"ahler hyperbolicity and $L\sb 2$-Hodge theory}, 
J. Diff. Geom. (1991), 263-292. 
\smallskip
\item{[Gr3]} M. Gromov, {\it Asymptotic invariants of infinite groups}, in {\it Geometric group theory II (Sussex, 1991)}, 1-295, Lond.
Math. Soc. Lect. Note 182, 1993.
\smallskip
\item{[GLP]} M. Gromov, J. Lafontaine, P. Pansu, {\it Structure m\'etriques pour les vari\'et\'es riemanniennes}, CEDIC/Nathan 1981.
\smallskip
\item{[L1]} U. Lang, {\it Higher-dimensional linear isoperimetric inequalities in hyperbolic groups}, Int. Math. Res. Notices 2000, 709-717.
\smallskip
\item{[M1]} I. Mineyev, {\it Higher dimensional isoperimetric functions in hyperbolic groups}, Math. Z. 233
(2000), 327-345.
\smallskip
\item{[M2]} I. Mineyev, {\it Bounded cohomology characterizes hyperbolic groups}, Quart. J. Math. 53 (2002), 59-73.
\smallskip
\item{[dR]} G. de Rham, {\it Vari\'et\'es diff\'erentiables}, Hermann, 1955 (English translation:
{\it Differentiable manifolds}, Springer GMW 266, 1984.
\smallskip
\item{[Si]} J.-C. Sikorav, {\it Growth of a primitive of a differential form}, Bull. Soc. Math. France 129 (2001), 159-168.
\smallskip
\item{[Su]} D. Sullivan, {\it Cycles for the dynamical study of foliated manifolds}, Invent. Math. 36 (1976), 225-255.
\smallskip
\item{[W]} H. Whitney, {\it Geometric integration theory}, Princeton Math Series 21, 1957.

\end
One can in fact find a primitive on $V$ which is almost optimal for all $R$, as follows.
If  $\psi:R_+\to\R_+^*$ be any continuous function such that $\psi(R) >{\rm Cof}(R)$ for all $R$, let $F:\R_+\to\R_+^*$ be such that 
${\rm Cof}<F<\psi$. 

Let $\omega$ be an exact  $2$-form $\omega$ with $||\omega||\le1$.
For every $R\ge0$, the definition of ${\rm Cof}(R)$ gives a primitive $\tau_R$ such that $||\tau_R||\le F(R)$ on $B'(x_0,R)$. Thus for every $T$ with support in $B(x_0,R)$, one has
$$|T(\omega)|=|\partial T(\tau_R)|\le F(R){\rm M}(\partial T)\le {\rm M}_\varphi(\partial T),$$
where $\varphi(x)=F(d(x_0,x))$. Theorem 2 implies that $\omega$ has a primitive such that $||\tau||(x)\le \psi(d(x_0,x))$ for every $x$.
{\bf N'importe quoi, en fait c'est une question ouverte !}

\medskip
(ii) 

\medskip
{\bf Corollary.}\hskip3mm {\it If all the relations in $R$ are triangular, $$\eqalign{\Delta^{ab}_{\hbox{\bbsmall R}}(w)&=\max\{t(I_w)\mid t\in C^1(G;\R),||u_{dt}||\le1\}\cr
&=\max\{\sum_{i=1}^n t(g_i,g_{i+1})\mid , |t(g\overline s_1,g\overline s_1\overline s_2)-t(g,g\overline s_1\overline s_2)+t(g,g\overline s_1)|\le1\}.\cr}$$
}

$$\eqalign{C^0(\widetilde K)=\R[G]^{c_0}\ , \ C_0(G)=\R[G]\cr
C_1(\widetilde K)=\R[G]^{c_1} \ , \ C_1(G)=\R[G^2]/\hbox{\got S}_2\cr
C_2(\widetilde K)=\R[G]^{c_2}\ ,\  C_2(G)=\R[G^3]/\hbox{\got S}_3\cr}$$

\medskip
{\bf Proof.}\hskip3mm This is nothing more than adding bounds to the proof of the theorem of de Rham given in [ST], p.165 sqq.
\medskip
(i) This is Lemma 1 of ([ST], p.165) , with bounds added (cf also [Si], p.163). The map $R^k$ is defined by
$R^k(u)=\sum_{\sigma\in K_q}\alpha_\sigma$, where 
$$\alpha_\sigma=k!\sum_{i=0}^qg_{j_i}dg_{j_0}\wedge\cdots\wedge\widehat{dg_{j_i}}\wedge dg_{j_k},$$
where $(g_j)$ is a partition of unity subordinated to ${\rm st}(v_j)$, with $||Dg_j||$ bounded.
\medskip
(ii) Let $U^k$ be a regular neighbourhood of the standard simplex $\Delta^k\subset\R^n$, and let $U^k_\partial\subset U$ be a regular neighbourhood of $\partial\Delta^k$, which is a collar of a part of $\partial U^k$ (not all of $\partial U^k$ as wrongly stated in [Si]!). 
\smallskip
Let ${\cal E}_k^q={\cal E}^k_0(U^k,U^k_\partial)$, $({\cal E}_k^q)_0=\{\alpha\in{\cal E}^q_k\mid\int_{\Delta^k}\alpha=0\}$, ${\cal Z}^q_k=\{\alpha\in{\cal E}^q_k\mid d\alpha=0\}$, ${\cal B}^q_k={\cal Z}^q_k$ if $k\ne q$, and ${\cal B}^q_q={\cal Z}^q_q\cap({\cal E}^q_q)_0$.
\smallskip
Note that $d$ sends ${\cal E}^q_k$ into
${\cal B}^q_k$.

\medskip

{\bf Lemma 1.}\hskip3mm {\it  There exists a primitivation operator $P^q_k:{\cal B}^q_k\to {\cal E}^{q-1}_k$ which is bounded in $L^\infty$ and has an image contained in $({\cal E}^{q-1}_{q-1})_0$ if $k=q-1$.}
\medskip
{\bf Proof of the lemma.}\hskip3mm This is essentially Lemma 3 in ([ST] p.169) with bounds. In [Si], this lemma is proved but for the property that the image of $R^q_{q-1}$ is contained in $({\cal E}^q_k)_0$. This last property is easy to achieve: fix $\tau_k\in{\cal E}^k_q$  such that $\int_{\Delta^k}\tau_k=1$, and if
$P'^q_{q-1}$ is a bounded right inverse, define
$$P^q_{q-1}\omega=P'^q_{q-1}\omega-\big(\int_{\Delta^{q-1}}P'^q_{q-1}\omega\big)\tau_{q-1}.$$
\medskip

For each simplex $s\in K$ there is a neighbourhood $U(s)$ and a bounded diffeomorphism $\varphi_s:U(s)\to U^k$, $k=\dim s$. Then $U(s)_\partial=\varphi_s^{-1}(U^k_\partial)$ is a neighbourhood of $\partial s$. We choose $(U^k)$ ``thinner and thinner'' such 
that $U(s)_\partial\subset\!\subset \bigcup_{t\subset\partial s} U(t)$ and $U(s)\cap U(s')\subset U(\partial s)$ if $t\subset s\cap s'$ (rather than $U(s)_\partial$ as in [Si])).
\medskip
Let
${\cal U}^k=\cup_{s\in K_k} U(s)$, which is a neighbourhood of the $k$-skeleton $K^{(k)}$, and ${\cal V}^k=\cup_{s\in K_k} U(s)_\partial$, which is a neighbourhood
of $K^{(k-1)}$ compactly contained in ${\cal U}_{k-1}$. We can  
find a smooth function $\rho_k:M\to [0,1]$ such that
 $\rho_k=1$ on ${\cal V}_k$, $\rho_k=0$ outside ${\cal U}_{k-1}$, and $||D\rho_k||$ is bounded.

\smallskip
Let ${\cal E}_k^q={\cal E}^q({\cal U}^k,{\cal V}^k)$, $({\cal E}_k^q)_0={\cal E}^q_k\cap\ker I^k$, ${\cal Z}^q_k=\{\alpha\in{\cal E}^q_k\mid d\alpha=0\}$, ${\cal B}^q_k={\cal Z}^q_k$ if $k\ne  q$, and ${\cal B}^q_q={\cal Z}^q_q\cap \ker I^q$.

\smallskip
Note that $d$ sends ${\cal E}^q(U^k,U^k_\partial)$ into
${\cal B}^q(U^k,U^k_\partial)$.

\medskip
{\bf Lemma 2.}\hskip3mm {\it There exists a primitivation operator ${\cal P}^q_k:{\cal B}_k^q\to{\cal E}_k^{q-1}$ such that
$$||P^q_k\alpha||_{L^\infty(U(s))}\le C||\alpha||_{L^\infty(U(s))}$$
and which has an image contained in $({\cal E}^{q-1}_{q-1})_0$ if $k=q-1$.}
\medskip
{\bf Proof.}\hskip3mm It suffices to glue together the operators :
$(P^q_k)_s=\varphi_s^*\circ P^q_k\circ\varphi_s^{-1*}$
 associated to the simplices by Lemma 1 and the diffeomorphisms $\varphi_s$.
 
 \medskip
 {\bf End of the proof.}\hskip3mm  We define inductively $\tau_0\in{\cal E}^{q-1}(V)$  by
 $\tau_{-1}=0$ and 
 $$\tau_k=\tau_{k-1}+\rho_k{\cal P}_k^q(\omega-d\tau_{k-1}).$$ This is well defined  
 since (by induction) $d\tau_{k-1}=\omega$ on ${\cal U}_{k-1}$ and a fortiori on ${\cal V}_k$, and (for $k=q$)
 $I^q(\omega)=0$, $(I^{q-1}(\tau_{q-1})=0)\Rightarrow I^q(d\tau^{q-1})=0)$. For $q=n$, $d\tau_n=\omega$ on $V$. Clearly, $\tau_n=\Pi(\omega)$ depends linearly on $\omega$.
 \smallskip
 Finally, the estimate of Lemma 2,  the bound on $||D\rho_k||$, and the fact that ${\rm diam}(U(s))$ is bounded imply inductively that, if $x\in U(s)$,
 $$||\tau_k(x)||+||d\tau_k(x)||\le C_k\max\{||\omega(y)||\mid y\in B'(x,C_k\}.$$
 For $k=n$, this finishes the proof of the proposition.

 We cut $\gamma$ in $N$ consecutive arcs $I_i=\gamma|[t_{i-1},t_i]$, as follows: 
$t_0=0$, 
$$t_i=t_{i-1}+{1\over 2}d(x_0,\gamma(t_{i-1}))$$ as long as it is possible, ie $t_{i+1}\le\ell(\gamma)$. Since $d(x_0,\gamma(t_i)\ge1$, there exists a maximal $i$ ($\le 2\ell(\gamma)$) such that $t_i+{1\over 2}d(x_0,\gamma(t_i))\le\ell(\gamma)$. We define $N=i+2$ and $t_N=\ell(\gamma)$. Note that by construction,
$$\forall i\in[[1,N-2]])\hskip5mm \ell(I_i)={1\over2} d(x_0,\gamma(t_{i-1})).$$
There remains to define $t_{N-1}$. Let $t'_{N-1}=t_{N-2}+{1\over2}d(x_0,\gamma(t_{N-2}))$, then  $t'_{N-1}+{1\over2}d(x_0,\gamma(t'_{N-1}))>\ell(\gamma)$, thus
$$t_N-t'_{N-1}<{1\over2}d(x_0,\gamma(t'_{N-1}))<{1\over4}d(x_0,\gamma(t_{N-2})).$$
\smallskip
1) If $t_N-t'_{N-1}\ge{1\over 4}d(x_0,\gamma(t'_{N-1}))$, we define
$t_{N-1}=t'_{N-1}$. Then 
$$\eqalign{\ell(I_{N-1})&={1\over 2}d(x_0,\gamma(t_{N-2}))\cr
 {\ell(I_N)\over d(x_0,\gamma(t_{N-2}))}&\in[{1\over4},{1\over2}[.\cr}$$

\medskip
2) If $t_N-t'_{N-1}<{1\over 4}d(x_0,\gamma(t'_{N-1}))$, we define
$t_{N-1}=t_{N-2}+{1\over 4}d(x_0,\gamma(t_{N-2}))$. Then
$$d(x_0,\gamma(t_{N-1}))\le{5\over4}d(x_0,\gamma(t_{N-2}))$$
and
$$\eqalign{\ell(I_{N-1})&={1\over4}d(x_0,\gamma(t_{N-2}))\cr
 \ell(I_N)&\ge{1\over4} d(x_0,\gamma(t_{N-2}))\ge{1\over5}d(x_0,\gamma(t_{N-1}))\cr
\ell(I_N)&\le(t_N-t'_{N-1})+{1\over4}d(x_0,\gamma(t_{N-2}))\cr&
\le{1\over 2}d(x_0,\gamma(t_{N-2})).\cr}$$
\medskip
Thus $\displaystyle{\ell(I_i)\over d(x_0,\gamma(t_{i-1}))}\in[{1\over5},{1\over2}]$ for all $i$. Let $c_i$ be a minimal  geodesic from $x_0$ to $\gamma(t_i)$, and let $\gamma_i=c-{i-1}*\gamma|I_i*c_i^{-1}$. Set $d_i=d(x_0,\gamma(t_i)), \ell_i=\ell(\gamma_i)=t_i-t_{i-1},\delta_i=d(x_0,\gamma(I_i))$. Then $\displaystyle{\ell_i\over d_{i-1}}\in[{1\over5},{1\over2}]$, thus
$$\eqalign{\delta_i&\ge d_{i-1}-\ell_i\ge{1\over2}d_{i-1}\ge2\ell_i\cr
d_i&\le d_{i-1}+\ell_i\le {3\over2}d_{i-1}\le \delta_i.\cr}$$

$$\eqalign\R FA(\gamma)&\le\sum_{i=1}^N{\rm FA}^{\hbox{\bbsmall R}}(\gamma_i)\cr
&\le\sum_{i=1}^N{\rm FA}^{\hbox{\bbsmall R}}(d(x_0,\gamma(t_{i-1}))+(t_i-t_{i-1})+d(x_0,\gamma(t_i)))\cr
&\le\sum_{i=1}^N{\rm FA}^{\hbox{\bbsmall R}}(d(x_0,\gamma(t_{i-1}))+(t_i-t_{i-1})+d(x_0,\gamma(t_i)))\cr}$$

\medskip
{\bf Proposition.}\hskip3mm {\it Let $\omega$ be an exact $q$-form with norm $\le1$. Then $\omega$ has a locally flat primitive $\tau_1$ such that $||\tau||(x)\le{\rm Cof}(d(x_0,x))$.

Thus, for every $g>{\rm Cof}$, $\omega$ has a smooth primitive $\tau$ such that $\||\tau(x)||\le g(d(x_0,x))$.}
\medskip
{\bf Proof.}\hskip3mm Let $\chi(x)={\rm Cof}(d(x_0,x))$, $\psi(x)=g(d(x_0,x))$, and $\varphi=(1/2)(\chi+\psi)$. By Theorem 1, it suffices to prove that $T(\omega)\le{\rm M}_\varphi(\partial T)$ for every normal $T$. 
\smallskip
Fix $R\ge0$. The definition of ${\rm Cof}$ implies that there exist a primitive $\tau_R$ which satisfies $||d\tau_R||\le {\rm Cof}(R)$ on $B(x_0,R)$. Thus

\medskip
We also give a localized version.
\medskip
{\bf Theorem 1'.}\hskip3mm {\it Let $\omega\in{\cal E}^q(V)$ with $q\ge2$, and let $\varphi:A\to\R_+$ be continuous, where $A\subset V$ is closed. Let $U$ be a any 
open neighborhood of $A$.
Assume that $I_c(\omega)\le{\rm M}_\varphi(I_{\partial c})$ for every real smooth singular $q$-chain $c$ on $V$ with boundary in $U$. 
\smallskip
Then for every continuous function $\psi>\varphi$ on $A$, there exists a smooth form $\tau\in{\cal E}^{q-1}(V)$ such that $d\tau=\omega$ and $||\tau||\le \psi$ on $A$.}
\medskip

\bigskip
{\bf Proof of Theorem 1'.}\hskip3mm Similarly, one has $T(\omega)\le {\rm M}_\varphi(\partial T)$
is satisfied for every $T\in {\cal E}'_q(V)$ with ${\rm supp}(\partial T)\subset U$.
Thus there exists a locally flat form $\tau_0$ such that  $d\tau_0=\omega$ and
$S(\tau_0)\le {\rm M}_\varphi(S)$ for every $S\in {\cal E}'_{q-1}(V)$ with ${\rm supp}(S)\subset U$. Define 
$\tau={\cal R}^*_\varepsilon\tau_0-{\cal H}^*_\varepsilon\omega=\tau_0+d{\cal H}^*_\varepsilon \tau_0.$
This proves Theorem 1'.

our stronger inequalities without saying so. 
\smallskip
Indeed, he proceeds as follows: By [F1],  Theorem 4.1.23, , there are chains $P_1,S_1,P_2,S_2$ with support  in $K$ such that $P_1,P_2$ are polyhedral and
$$ \eqalign{F_K(P_1-T)&<\varepsilon/4\ , \ {\rm M}(P_1)<{\rm M}(T)+\varepsilon/4\cr
F_K(P_2-\partial T)&<\varepsilon/4\ , \ {\rm M}(P_2)<{\rm M}(\partial T)+\varepsilon/4.\cr}$$
Then $F_K(P_2-\partial P_1)<3\varepsilon/4$. 
By [F1], Lemma 4.2.23, there is a polyhedral chain $Y$ with support in $K$, such that 
$${\rm M}(P_2-\partial P_1-\partial Y)+{\rm M}(Y)<3\varepsilon/4.$$ Then $P=P_1+Y$ is a polyhedral chain with support in $K$, and satisfies
$$\eqalign{{\rm M}(P)&\le {\rm M}(P_1)+{\rm M}(Y)< {\rm M}(T)+\varepsilon\cr
{\rm M}(\partial T')&\le{\rm M}(\partial P_1+\partial Y-P_2)+{\rm M}(P_2)\le {\rm M}(P_1)+{\rm M}(Y)< {\rm M}(\partial P)+\varepsilon.\cr
}$$
Finally, $P-T-\partial S_1=(P_1-T-\partial S_1)+Y$ which has a mass $<\varepsilon/2$, finishing the proof for $V=\R^n$.

\medskip
{\bf Proof.}\hskip3mm We follow the book of de Rham [dR] (p.72-83), adding the pointwise estimates.
First, we regularize flat forms in the Euclidean space. Fix a smooth positive function $\rho$ on $\R^n$, with support in $B(1)$, such that $\int_{\R^n}\rho(x)dx=1$.
For $\omega\in{\cal E}^q(\R^n)$, we define, for every $\varepsilon>0$:
$$\eqalign{R_\varepsilon^*\omega(x)&=\int_{\R^n}\varepsilon^{-n}\rho(\varepsilon^{-1}y)\omega(x+y)\ dy\cr
H_\varepsilon^*\omega(x)&=\int_{\R^n}\varepsilon^{-n}\rho(\varepsilon^{-1}y)\ (\int_0^1\iota_y\omega(x+ty)\ dt)\ dy.\cr}$$

Then $R^*_\varepsilon\omega\in{\cal E}^k(\R^n)$, $H^*_\varepsilon\omega\in{\cal E}^{k-1}(\R^n)$. If $\omega\in{\cal E}^*(\R^n)$ $\varepsilon\to0$, ${\cal R}^*_\varepsilon\omega\to\omega$ in $C^{\infty}_{loc}$.

If $\omega$ is smooth, easy computations give
$$\leqalignno{R^*_\varepsilon\omega-\omega&=d(H^*_\varepsilon\omega)+H^*_\varepsilon d\omega&(1)\cr
d(R^*_\varepsilon\omega)&=R^*_\varepsilon(d\omega)=d\omega+H^*_\varepsilon d(H^*_\varepsilon\omega).&(2)\cr}$$
Assume now that $\omega$ is flat, ie $d\omega\in L^\infty$ in the sense of distributions, or
$|\int_{\R^n}\omega\wedge d\varphi|\le C||\varphi|L_{L^\infty}$. Then integrations by parts show that (1) and (2) still hold in the sense of distributions (cf [dR] p.75). For instance, the first identity means $\int_{\R^n}(R^*_\varepsilon\omega-\omega)\wedge \varphi=\int_{\R^n}(-1)^kH^*_\varepsilon\omega\wedge d\varphi+H^*_\varepsilon d\omega\wedge \varphi$.
\smallskip
Furthermore, for the Euclidean metric one has
$$\leqalignno{||R^*_\varepsilon\omega(x)||&\le||\omega|B(x,\varepsilon)||&(3)\cr
||H^*_\varepsilon\omega(x)||&\le \varepsilon\ ||\omega|B(x,\varepsilon)||.&(4)\cr}$$

\medskip
We now turn to the Riemannian manifold $V$. Let $(h_i)$ be a smooth partition of unity on $V$ such that $h_i=1$ on $U_i$ and $h_i=0$ outside of $V_i$. 

Let $(\varepsilon_i)$ be a sequence of positive numbers, which we shall later require to decrease sufficiently fast. 
Let
${\cal R}^{*i}_\varepsilon=R^*_{i,\varepsilon_i}\cdots R^*_{0,\varepsilon_0}$ and
${\cal H}^i_\varepsilon=H^{*i}_{i,\varepsilon_i}R^*_{i-1,\varepsilon_{i-1}}\cdots R^*_{0,\varepsilon_0}.$
For $\omega\in{\cal F}^k_{loc}(V)$ and $\varepsilon_i>0$, define
$$\eqalign{R^*_{i,\varepsilon_i}\omega&=\varphi_i^*R^*_\varepsilon\varphi_i^{-1*}(h_i\omega)+(1-h_i)\omega\cr
H^*_{i,\varepsilon_i}\varphi&=\varphi_i^*H^*_\varepsilon\varphi_i^{-1*}(h_i\omega).\cr}$$
Then $R^*_{i,\varepsilon_i}\omega$ is a locally flat form, which is smooth on $U_i$ and coincides with $\omega$ outside
$V_i$. The form $H^*_{i,\varepsilon_i}\omega$ is smooth everywhere.
\smallskip
Appplying (1) and (2) to $\varphi_i^{-1*}(h_i\omega)$, we get
$$\leqalignno{R^*_{i,\varepsilon_i}\omega-\omega&=d(H^*_{i,\varepsilon_i}\omega)+H^*_{i,\varepsilon_i}d\omega&(5)\cr
R^*_{i,\varepsilon_i}(d\omega)&=d(R^*_{i,\varepsilon_i}\omega)&(6)\cr}$$
Using (3), (4) and $||d\varphi_i^{\pm1}||\le C_i$, we have
$$\leqalignno{||R^*_{i,\varepsilon}\omega(x)||&\le(1+C_i^2\varepsilon)||h_i\omega|B(x,C_i\varepsilon))||+(1-h_i(x))||\omega|B(x,C_i\varepsilon)||\cr
&\le(1+C'_i\varepsilon)\ ||\omega|B(x,C_i\varepsilon)||&(7)\cr
||H^*_{i,\varepsilon}\omega(x)||&\le C_i^2\varepsilon\ ||\omega|B(x,C_i\varepsilon)||.&(8)\cr}$$

Outside $V_i$, $R^*_{i,\varepsilon_i}$ is the identity and $H^*_{i,\varepsilon_i}$ vanishes, thus
$$\matrix{{\cal R}^*_\varepsilon=\lim_i R^{*i}_\varepsilon&{\rm and}&{\cal H}_\varepsilon=\sum_i H^{*i}_\varepsilon\cr}$$
are well defined, and for every compact $K\subset V$ there exists $N={\rm N}(K)\in\N$ such that
${\cal R}^*_\varepsilon|K=R^N_\varepsilon$ and ${\cal H}^*_\varepsilon|K=\sum_{i=0}^N H^i_\varepsilon$.
Thus ${\cal R}^*_\varepsilon\omega$ is smooth everywhere on $V$, and ${\cal R}^*_\varepsilon$ commutes with $d$ by (6).
\smallskip
By (5), 
$R^{*i}_\varepsilon\omega-R^{*i-1}_\varepsilon\omega=dH^{*i}_\varepsilon\omega+H^{*i}_\varepsilon d\omega.$
Adding this for all $i$ we get
${\cal R}^*_\varepsilon\omega-\omega=d{\cal H}^*_\varepsilon\omega+{\cal H}^*_\varepsilon d\omega.$
\smallskip
For $x\in V$, denote by $I(x)$ the (finite and locally bounded) set of $i$ such that $x\in V_i$. The estimates (7) and (8) give
$$\eqalign{||{\cal R}^*_\varepsilon\omega(x)||&\le\big(\prod_{i\in I(x)}(1+C_i\varepsilon_i)\big)\ ||\omega|B(x,\sum_{i\in I(x)}\varepsilon_i)|||_{L^\infty}\cr
||{\cal H}^*_\varepsilon\omega(x)||&\le\big(\sum_{i\in I(x)}C_i\varepsilon_i\big)\ ||\omega|B(x,\sum_{i\in I(x)}\varepsilon_i))||_{L^\infty}.\cr}$$
If $(\varepsilon_i)$ decreases sufficiently fast, this gives the announced inequalities \smallskip
The last statement of Proposition 2 is immediate since only a finite number of $V_i$ are involved and the associated $\varepsilon_i$ tend to zero.

\medskip
{\bf Proof.}\hskip3mm Use Proposition 2 with $\varepsilon$ constant, and define ${\cal R}_\varepsilon T(\omega)=T({\cal R}^*_\varepsilon\omega)$, ${\cal H}_\varepsilon T(\omega)=T({\cal H}^*_\varepsilon\omega).$
The algebraic properties are immediate. The weak limit and the upper bound on ${\rm M}({\cal R}_\varepsilon T)$ result from 
$\lim_{\varepsilon\to0}{\cal R}^*_\varepsilon\omega=\omega
$.

This implies that $t$ is a locally flat cochain, thus there exists a locally flat form $\tau_0\in {\cal F}^{q-1}(V)$ such that $t(S)=S(\tau_0)$
if ${\rm M}(S)<\infty$. 
Thus there exists a locally flat form $\tau_0$ such that  $d\tau_0=\omega$ and
$S(\tau_0)\le {\rm M}_\varphi(S)$ for every $S\in {\cal E}'_{q-1}(V)$ with ${\rm supp}(S)\subset U$. Define 
$\tau={\cal R}^*_\varepsilon\tau_0-{\cal H}^*_\varepsilon\omega=\tau_0+d{\cal H}^*_\varepsilon \tau_0.$
This proves Theorem 1'.

\medskip
{\bf Proof of Theorem 1.}\hskip3mm By the density theorem for singular chains, $T(\omega)\le {\rm M}_\varphi(\partial T)$ for every $T\in {\cal E}'_q(V)$.
By Hahn-Banach, there exists a linear form $t$ on ${\cal E}'_{q-1}(V)$ such that $t(\partial T)=T(\omega)$ for every $T\in {\cal E}'_q(V)$ and 
$t(S)\le {\rm M}_\varphi(S)$ for every $S\in {\cal E}'_{q-1}(V)$. Thus
$$t(S)=t(S-\partial T)+t(\partial T)\le {\rm M}_\varphi(S-\partial T)+M_{||\omega||}(T).$$
Using the regularization procedure of section 1, define 
$$\tau={\cal R}^*_\varepsilon\tau_0-{\cal H}^*_\varepsilon\omega=\tau_0+d{\cal H}^*_\varepsilon \tau_0.$$
Then $\tau$ is smooth and $d\tau=d\tau_0=\omega$. Moreover, if $x\in U$,
$$||\tau(x)||\le(1+\varepsilon(x))\big)\ ||\varphi|B(x,\varepsilon(x))||_{L^\infty}+
\varepsilon(x)\ ||\omega|B(x,\varepsilon(x))||_{L^\infty}.$$
If $\psi>\varphi$ is a given continuous function, for $\varepsilon$ decreasing sufficiently fast the right-hand side above
is $\le \psi(x)$, qed.

\medskip
{\bf Proposition 2 (regularization of currents).}\hskip3mm {\it For every $\varepsilon>0$, there exists a linear chain map of degree $0$, ${\cal R}_\varepsilon:{\cal E}'(V)\to{\cal S}'(V)$, and a homotopy 
${\cal R}_\varepsilon-{\rm Id}=\partial{\cal H}_\varepsilon+{\cal H}_\varepsilon \partial$,
such that if $T$ has finite mass, one has the weak limit
$\displaystyle\lim_{\varepsilon\to0}{\cal R}_\varepsilon T=T$, and
$\displaystyle\limsup_{\varepsilon\to0}{\rm M}({\cal R}_\varepsilon T)={\rm M}(T).$
}